\documentclass[12pt]{amsart}
\usepackage{amssymb}
\input{xypic}
\begin{document}
\theoremstyle{plain}
\newtheorem{thm}{Theorem}[section]
\newtheorem*{thm1}{Theorem 1}
\newtheorem*{thm2}{Theorem 2}
\newtheorem{lemma}[thm]{Lemma}
\newtheorem{lem}[thm]{Lemma}
\newtheorem{cor}[thm]{Corollary}
\newtheorem{prop}[thm]{Proposition}
\newtheorem{propose}[thm]{Proposition}
\newtheorem{variant}[thm]{Variant}
\theoremstyle{definition}
\newtheorem{notations}[thm]{Notations}
\newtheorem{rem}[thm]{Remark}
\newtheorem{rmk}[thm]{Remark}
\newtheorem{rmks}[thm]{Remarks}
\newtheorem{defn}[thm]{Definition}
\newtheorem{ex}[thm]{Example}
\newtheorem{claim}[thm]{Claim}
\newtheorem{ass}[thm]{Assumption}
\numberwithin{equation}{section}
\newcounter{elno}                
\def\points{\list
{\hss\llap{\upshape{(\roman{elno})}}}{\usecounter{elno}}} 
\let\endpoints=\endlist


\catcode`\@=11
%
%
\def\opn#1#2{\def#1{\mathop{\kern0pt\fam0#2}\nolimits}} 
\def\bold#1{{\bf #1}}%
\def\underrightarrow{\mathpalette\underrightarrow@}
\def\underrightarrow@#1#2{\vtop{\ialign{$##$\cr
 \hfil#1#2\hfil\cr\noalign{\nointerlineskip}%
 #1{-}\mkern-6mu\cleaders\hbox{$#1\mkern-2mu{-}\mkern-2mu$}\hfill
 \mkern-6mu{\to}\cr}}}
\let\underarrow\underrightarrow
\def\underleftarrow{\mathpalette\underleftarrow@}
\def\underleftarrow@#1#2{\vtop{\ialign{$##$\cr
 \hfil#1#2\hfil\cr\noalign{\nointerlineskip}#1{\leftarrow}\mkern-6mu
 \cleaders\hbox{$#1\mkern-2mu{-}\mkern-2mu$}\hfill
 \mkern-6mu{-}\cr}}}
%
%

%
\def\:{\colon}
\let\oldtilde=\tilde
\def\tilde#1{\mathchoice{\widetilde{#1}}{\widetilde{#1}}%
{\indextil{#1}}{\oldtilde{#1}}}
\def\indextil#1{\lower2pt\hbox{$\textstyle{\oldtilde{\raise2pt%
\hbox{$\scriptstyle{#1}$}}}$}}
\def\pnt{{\raise1.1pt\hbox{$\textstyle.$}}}
%

%
\let\amp@rs@nd@\relax
\newdimen\ex@
\ex@.2326ex
\newdimen\bigaw@
\newdimen\minaw@
\minaw@16.08739\ex@
\newdimen\minCDaw@
\minCDaw@2.5pc
\newif\ifCD@
\def\minCDarrowwidth#1{\minCDaw@#1}
\newenvironment{CD}{\@CD}{\@endCD}
\def\@CD{\def\A##1A##2A{\llap{$\vcenter{\hbox
 {$\scriptstyle##1$}}$}\Big\uparrow\rlap{$\vcenter{\hbox{%
$\scriptstyle##2$}}$}&&}%
\def\V##1V##2V{\llap{$\vcenter{\hbox
 {$\scriptstyle##1$}}$}\Big\downarrow\rlap{$\vcenter{\hbox{%
$\scriptstyle##2$}}$}&&}%
\def\={&\hskip.5em\mathrel
 {\vbox{\hrule width\minCDaw@\vskip3\ex@\hrule width
 \minCDaw@}}\hskip.5em&}%
\def\verteq{\Big\Vert&&}%
\def\noarr{&&}%
\def\vspace##1{\noalign{\vskip##1\relax}}\relax\let\amp@rs@nd@&\iffalse}\fi
 \CD@true\vcenter\bgroup\relax\let\\=\cr\iffalse}\fi\tabskip\z@skip\baselineskip20\ex@
 \lineskip3\ex@\lineskiplimit3\ex@\halign\bgroup
 &\hfill$\m@th##$\hfill\cr}
\def\@endCD{\cr\egroup\egroup}
%
\def\>#1>#2>{\amp@rs@nd@\setbox\z@\hbox{$\scriptstyle
 \;{#1}\;\;$}\setbox\@ne\hbox{$\scriptstyle\;{#2}\;\;$}\setbox\tw@
 \hbox{$#2$}\ifCD@
 \global\bigaw@\minCDaw@\else\global\bigaw@\minaw@\fi
 \ifdim\wd\z@>\bigaw@\global\bigaw@\wd\z@\fi
 \ifdim\wd\@ne>\bigaw@\global\bigaw@\wd\@ne\fi
 \ifCD@\hskip.5em\fi
 \ifdim\wd\tw@>\z@
 \mathrel{\mathop{\hbox to\bigaw@{\rightarrowfill}}\limits^{#1}_{#2}}\else
 \mathrel{\mathop{\hbox to\bigaw@{\rightarrowfill}}\limits^{#1}}\fi
 \ifCD@\hskip.5em\fi\amp@rs@nd@}
\def\<#1<#2<{\amp@rs@nd@\setbox\z@\hbox{$\scriptstyle
 \;\;{#1}\;$}\setbox\@ne\hbox{$\scriptstyle\;\;{#2}\;$}\setbox\tw@
 \hbox{$#2$}\ifCD@
 \global\bigaw@\minCDaw@\else\global\bigaw@\minaw@\fi
 \ifdim\wd\z@>\bigaw@\global\bigaw@\wd\z@\fi
 \ifdim\wd\@ne>\bigaw@\global\bigaw@\wd\@ne\fi
 \ifCD@\hskip.5em\fi
 \ifdim\wd\tw@>\z@
 \mathrel{\mathop{\hbox to\bigaw@{\leftarrowfill}}\limits^{#1}_{#2}}\else
 \mathrel{\mathop{\hbox to\bigaw@{\leftarrowfill}}\limits^{#1}}\fi
 \ifCD@\hskip.5em\fi\amp@rs@nd@}
%
%
\newenvironment{CDS}{\@CDS}{\@endCDS}
\def\@CDS{\def\A##1A##2A{\llap{$\vcenter{\hbox
 {$\scriptstyle##1$}}$}\Big\uparrow\rlap{$\vcenter{\hbox{%
$\scriptstyle##2$}}$}&}%
\def\V##1V##2V{\llap{$\vcenter{\hbox
 {$\scriptstyle##1$}}$}\Big\downarrow\rlap{$\vcenter{\hbox{%
$\scriptstyle##2$}}$}&}%
\def\={&\hskip.5em\mathrel
 {\vbox{\hrule width\minCDaw@\vskip3\ex@\hrule width
 \minCDaw@}}\hskip.5em&}
\def\verteq{\Big\Vert&}
\def\novarr{&}
\def\noharr{&&}
\def\SE##1E##2E{\slantedarrow(0,18)(4,-3){##1}{##2}&}
\def\SW##1W##2W{\slantedarrow(24,18)(-4,-3){##1}{##2}&}
\def\NE##1E##2E{\slantedarrow(0,0)(4,3){##1}{##2}&}
\def\NW##1W##2W{\slantedarrow(24,0)(-4,3){##1}{##2}&}
\def\slantedarrow(##1)(##2)##3##4{%
\thinlines\unitlength1pt\lower 6.5pt\hbox{\begin{picture}(24,18)%
\put(##1){\vector(##2){24}}%
\put(0,8){$\scriptstyle##3$}%
\put(20,8){$\scriptstyle##4$}%
\end{picture}}}
\def\vspace##1{\noalign{\vskip##1\relax}}\relax\let\amp@rs@nd@&\iffalse}\fi
 \CD@true\vcenter\bgroup\relax\let\\=\cr\iffalse}\fi\tabskip\z@skip\baselineskip20\ex@
 \lineskip3\ex@\lineskiplimit3\ex@\halign\bgroup
 &\hfill$\m@th##$\hfill\cr}
\def\@endCDS{\cr\egroup\egroup}
%
\newdimen\TriCDarrw@
\newif\ifTriV@
\newenvironment{TriCDV}{\@TriCDV}{\@endTriCD}
\newenvironment{TriCDA}{\@TriCDA}{\@endTriCD}
\def\@TriCDV{\TriV@true\def\TriCDpos@{6}\@TriCD}
\def\@TriCDA{\TriV@false\def\TriCDpos@{10}\@TriCD}
\def\@TriCD#1#2#3#4#5#6{%
\setbox0\hbox{$\ifTriV@#6\else#1\fi$}
\TriCDarrw@=\wd0 \advance\TriCDarrw@ 24pt
\advance\TriCDarrw@ -1em
\def\SE##1E##2E{\slantedarrow(0,18)(2,-3){##1}{##2}&}
\def\SW##1W##2W{\slantedarrow(12,18)(-2,-3){##1}{##2}&}
\def\NE##1E##2E{\slantedarrow(0,0)(2,3){##1}{##2}&}
\def\NW##1W##2W{\slantedarrow(12,0)(-2,3){##1}{##2}&}
\def\slantedarrow(##1)(##2)##3##4{\thinlines\unitlength1pt
\lower 6.5pt\hbox{\begin{picture}(12,18)%
\put(##1){\vector(##2){12}}%
\put(-4,\TriCDpos@){$\scriptstyle##3$}%
\put(12,\TriCDpos@){$\scriptstyle##4$}%
\end{picture}}}
\def\={\mathrel {\vbox{\hrule
   width\TriCDarrw@\vskip3\ex@\hrule width
   \TriCDarrw@}}}
\def\>##1>>{\setbox\z@\hbox{$\scriptstyle
 \;{##1}\;\;$}\global\bigaw@\TriCDarrw@
 \ifdim\wd\z@>\bigaw@\global\bigaw@\wd\z@\fi
 \hskip.5em
 \mathrel{\mathop{\hbox to \TriCDarrw@
{\rightarrowfill}}\limits^{##1}}
 \hskip.5em}
\def\<##1<<{\setbox\z@\hbox{$\scriptstyle
 \;{##1}\;\;$}\global\bigaw@\TriCDarrw@
 \ifdim\wd\z@>\bigaw@\global\bigaw@\wd\z@\fi
 \mathrel{\mathop{\hbox to\bigaw@{\leftarrowfill}}\limits^{##1}}
 }
 \CD@true\vcenter\bgroup\relax\let\\=\cr\iffalse}\fi
 \tabskip\z@skip\baselineskip20\ex@
 \lineskip3\ex@\lineskiplimit3\ex@
 \ifTriV@
 \halign\bgroup
 &\hfill$\m@th##$\hfill\cr
#1&\multispan3\hfill$#2$\hfill&#3\\
&#4&#5\\
&&#6\cr\egroup%
\else
 \halign\bgroup
 &\hfill$\m@th##$\hfill\cr
&&#1\\%
&#2&#3\\
#4&\multispan3\hfill$#5$\hfill&#6\cr\egroup
\fi}
\def\@endTriCD{\egroup}
\newcommand{\mc}{\mathcal}
\newcommand{\mb}{\mathbb}
\newcommand{\surj}{\twoheadrightarrow}
\newcommand{\inj}{\hookrightarrow}
\newcommand{\zar}{{\rm zar}}
\newcommand{\an}{{\rm an}} 
\newcommand{\red}{{\rm red}}
\newcommand{\codim}{{\rm codim}}
\newcommand{\rank}{{\rm rank}}
\newcommand{\Ker}{{\rm Ker \ }}
\newcommand{\Pic}{{\rm Pic}}
\newcommand{\Div}{{\rm Div}}
\newcommand{\Hom}{{\rm Hom}}
\newcommand{\im}{{\rm im}}
\newcommand{\Spec}{{\rm Spec \,}}
\newcommand{\Sing}{{\rm Sing}}
\newcommand{\sing}{{\rm sing}}
\newcommand{\reg}{{\rm reg}}
\newcommand{\Char}{{\rm char}}
\newcommand{\Tr}{{\rm Tr}}
\newcommand{\Gal}{{\rm Gal}}
\newcommand{\Min}{{\rm Min \ }}
\newcommand{\Max}{{\rm Max \ }}
\newcommand{\Alb}{{\rm Alb}\,}
\newcommand{\GL}{{\rm GL}\,}        
\newcommand{\Ext}{{\rm Ext}\,}
\newcommand{\ie}{{\it i.e.\/},\ }
\newcommand{\niso}{\not\cong}
\newcommand{\nin}{\not\in}
\newcommand{\soplus}[1]{\stackrel{#1}{\oplus}}
\newcommand{\by}[1]{\stackrel{#1}{\rightarrow}}
\newcommand{\longby}[1]{\stackrel{#1}{\longrightarrow}}
\newcommand{\vlongby}[1]{\stackrel{#1}{\mbox{\large{$\longrightarrow$}}}}
\newcommand{\ldownarrow}{\mbox{\Large{\Large{$\downarrow$}}}}
\newcommand{\lsearrow}{\mbox{\Large{$\searrow$}}}
\renewcommand{\d}{\stackrel{\mbox{\scriptsize{$\bullet$}}}{}}
\newcommand{\dlog}{{\rm dlog}\,}    
\newcommand{\longto}{\longrightarrow}
\newcommand{\vlongto}{\mbox{{\Large{$\longto$}}}}
\newcommand{\limdir}[1]{{\displaystyle{\mathop{\rm lim}_{\buildrel\longrightarrow\over{#1}}}}\,}
\newcommand{\liminv}[1]{{\displaystyle{\mathop{\rm lim}_{\buildrel\longleftarrow\over{#1}}}}\,}
\newcommand{\norm}[1]{\mbox{$\parallel{#1}\parallel$}}
\newcommand{\boxtensor}{{\Box\kern-9.03pt\raise1.42pt\hbox{$\times$}}}
\newcommand{\into}{\hookrightarrow}
\newcommand{\image}{{\rm image}\,}
\newcommand{\Lie}{{\rm Lie}\,}      
\newcommand{\CM}{\rm CM}
\newcommand{\sext}{\mbox{${\mathcal E}xt\,$}}  
\newcommand{\shom}{\mbox{${\mathcal H}om\,$}}  
\newcommand{\coker}{{\rm coker}\,}  
\newcommand{\sm}{{\rm sm}}
\newcommand{\tensor}{\otimes}
\renewcommand{\iff}{\mbox{ $\Longleftrightarrow$ }}
\newcommand{\supp}{{\rm supp}\,}
\newcommand{\ext}[1]{\stackrel{#1}{\wedge}}
\newcommand{\onto}{\mbox{$\,\>>>\hspace{-.5cm}\to\hspace{.15cm}$}}
\newcommand{\propsubset}
{\mbox{$\textstyle{
\subseteq_{\kern-5pt\raise-1pt\hbox{\mbox{\tiny{$/$}}}}}$}}
\newcommand{\veq}{\mbox{\large $\parallel$}}
\newcommand{\sA}{{\mathcal A}}
\newcommand{\sB}{{\mathcal B}}
\newcommand{\sD}{{\mathcal D}}
\newcommand{\sE}{{\mathcal E}}
\newcommand{\sF}{{\mathcal F}}
\newcommand{\sG}{{\mathcal G}}
\newcommand{\sH}{{\mathcal H}}
\newcommand{\sI}{{\mathcal I}}
\newcommand{\sJ}{{\mathcal J}}
\newcommand{\sK}{{\mathcal K}}
\newcommand{\sL}{{\mathcal L}}
\newcommand{\sM}{{\mathcal M}}
\newcommand{\sN}{{\mathcal N}}
\newcommand{\sO}{{\mathcal O}}
\newcommand{\sQ}{{\mathcal Q}}
\newcommand{\sR}{{\mathcal R}}
\newcommand{\sS}{{\mathcal S}}
\newcommand{\sT}{{\mathcal T}}
\newcommand{\sU}{{\mathcal U}}
\newcommand{\sV}{{\mathcal V}}
\newcommand{\sW}{{\mathcal W}}
\newcommand{\sX}{{\mathcal X}}
\newcommand{\sY}{{\mathcal Y}}

\newcommand{\cal}{\mathcal}
\newcommand{\cA}{\sA}
\newcommand{\cC}{{\mathcal C}}
\newcommand{\cF}{\sF}
\newcommand{\cG}{\sG}
\newcommand{\ccH}{\sH}
\newcommand{\cHH}{\sH}
\newcommand{\cI}{\sI}
\newcommand{\ccL}{\sL}
\newcommand{\cO}{\sO}
\newcommand{\cV}{\sV}
\newcommand{\cM}{\sM}

\newcommand{\A}{{\mathbb A}}
\newcommand{\B}{{\mathbb B}}
\newcommand{\C}{{\mathbb C}}
\newcommand{\D}{{\mathbb D}}
\newcommand{\E}{{\mathbb E}}
\newcommand{\F}{{\mathbb F}}
\newcommand{\G}{{\mathbb G}}
\newcommand{\HH}{{\mathbb H}}
\newcommand{\J}{{\mathbb J}}
\newcommand{\M}{{\mathbb M}}
\newcommand{\N}{{\mathbb N}}
\renewcommand{\P}{{\mathbb P}}
\newcommand{\Q}{{\mathbb Q}}
\newcommand{\R}{{\mathbb R}}
\newcommand{\T}{{\mathbb T}}
\newcommand{\U}{{\mathbb U}}
\newcommand{\V}{{\mathbb V}}
\newcommand{\W}{{\mathbb W}}
\newcommand{\X}{{\mathbb X}}
\newcommand{\Y}{{\mathbb Y}}
\newcommand{\Z}{{\mathbb Z}}
\newcommand{\sC}{\C}
\newcommand{\sP}{\P}
\newcommand{\sZ}{\Z}
\title[A Lefschetz (1,1) Theorem]  
{A Lefschetz (1,1)-Theorem for normal projective complex varieties}  
\author{J. Biswas} 
\address{Institute of Mathematical Sciences, Taramani, Chennai-600 113,
India}
\email{jishnu@imsc.ernet.in}
\author{V. Srinivas} 
\address{School of Mathematics, TIFR,
Homi Bhabha Road, Bombay-400005, India}
\email{srinivas@math.tifr.res.in}
\date{}
\maketitle

\section{Introduction}

Let $X$ be a projective variety over $\C.$ Let $X_{an}$ be the
analytic space associated to $X.$ Let $c_1 : Pic(X) \to H^{2}(X_{an},
\Z)$ be the map which associates to a line bundle (or equivalently
a Cartier divisor) on $X$ its cohomology class. We may identify the
N\'eron-Severi group $NS(X)$ with the image of $Pic(X)$ in
$H^{2}(X_{an},\Z)$ under the above map. 

If $X$ is smooth, then by the Hodge decomposition theorem, we know that 
$$H^{2}(X_{an}, \C) = H^{2,0}(X_{an}) \oplus H^{1,1}(X_{an}) \oplus
H^{0,2}(X_{an}).$$ 
Let $F^{1}H^{2}(X_{an},\C) = H^{2,0}(X_{an}) \oplus H^{1,1}(X_{an}).$  
The {\em Lefschetz theorem on $(1,1)$ classes} ([GH], [L]) states that if
$X$ is a smooth, projective variety, then 
$$NS(X) = \{\alpha \in H^{2}(X_{an},\Z) | \alpha_{\sC} \in
F^{1}H^{2}(X_{an},\C)\}.$$ 

If $X$ is an arbitrary singular variety then by [D],~Theorem~8.2.2 the  
cohomology groups of $X$ with $\Z$-coefficients carry mixed Hodge
structures. Hence it makes sense to talk of $F^{1}H^{2}(X_{an},\C)$ for
such a variety $X$. Spencer Bloch, in a letter to Jannsen [J, appendix A],
asks whether the  ``obvious'' extension of the Lefschetz $(1,1)$ 
theorem is true for singular projective varieties, i.e., is it true that
$$NS(X) = \{ \alpha \in H^{2}(X_{an}, \Z) | \alpha_{\sC} \in
F^{1}H^{2}(X_{an},\C)\} ?$$

Barbieri-Viale and Srinivas [BS1] gave a counterexample to this question.
Let $X$ be a surface defined by the homogenous equation $w(x^3 - y^2z) +
f(x,y,z) = 0$ in $\P^3_{\sC},$ where $x,y,z,w$ are homogenous coordinates
in $\P^{3}_{\sC}$ and $f$ is a ``general'' homogenous polynomial over $\C$
of degree $4.$ They showed that for such an $X$, 
\[NS(X)\propsubset\{\alpha \in H^{2}(X_{an},\Z) |
\alpha_{\sC} \in F^{1}H^{2}(X_{an}, \C)\}.\] 
In the same paper [BS1] the authors ask the following question.
Let $X$ be a complete variety over $\C$. Let $H^{1}(X, {\cal H}^{1}_{X})$
be the subgroup of $H^{2}(X_{an},\Z)$  consisting of Zariski-locally
trivial cohomology classes, i.e., $\eta \in H^{2}(X_{an},\Z)$ lies in
$H^{1}(X,{\cal H}^{1}_{X})$ if and only if there exists a finite open
cover $\{U_{i}\}$ of $X$ by Zariski open sets such that $\eta \mapsto 0$
under the restriction maps $H^{2}(X_{an}, \Z) \to H^{2}((U_i)_{an}, \Z)$
for all $i$. Is
$$NS(X) = \{\alpha \in H^{1}(X,{\cal H}_{X}^{1}) | \alpha_{\sC} \in
F^{1}H^{2}(X_{an},\C)\} ?$$ 

We remark that for a smooth, projective variety $X$, if a cohomology class
$\eta \in H^{2}(X_{an},\Z)$ is zero when restricted to a nonempty Zariski
open set $U \subset X,$ then $\eta$ is the class of a divisor. So
$H^{1}(X,{\cal H}^{1}_{X}) = NS(X)$ and the above question has a positive
answer for a smooth, projective variety $X$. For any
projective variety $X$, there is an inclusion $NS(X)\subset H^{1}(X,{\cal
H}_{X}^{1})$;  Barbieri-Viale and Srinivas also give an example in [BS1]
of a singular variety for which this inclusion is strict. 

In general, for any projective variety $X$ over $\C,$
$$NS(X) \subset \{\alpha \in H^{1}(X,{\cal H}^{1}_{X}) | \alpha_{\sC}
\in F^1 H^2 (X_{an},\C)\}.$$ 
This follows from the inclusion $NS(X)\subset H^{1}(X,{\cal
H}_{X}^{1})$, combined with 
$$Ker(H^{2}(X_{an},\C) \to H^{2}(X_{an}, {\cal O}_{X_{an}})) \subset F^1
H^2 (X_{an},\C),$$ 
which is a consequence of results of Du Bois [DB] (and is also implicit
in [D]). When $X$ is normal, we prove the reverse inclusion, thereby
answering the question in the affirmative, for the normal case.
The statement of our Main Theorem is
\begin{thm}\label{(1,1)} Let $X$ be a normal, projective variety over
$\C.$ Then
$$NS(X) = \{ \alpha \in H^{1}(X, {\cal H}^{1}_{X}) | \alpha_{\sC} \in F^1
H^2 (X_{an},\C)\}.$$
\end{thm}
We also describe a counterexample, of a non-normal irreducible projective
3-fold with smooth normalization (isomorphic to
$\P^1\times\P^1\times\P^1$) for which the question has a {\em negative}
answer. However, it seems likely that the conclusion of the Theorem holds
for any {\it semi-normal} projective variety $X$ over $\C$. 

Our proof of Theorem~\ref{(1,1)} is in two steps: first we show that 
$H^{1}(X, {\cal H}^{1}_{X})\subset H^2(X_{an},\Z)$ is a sub-MHS of
level~1; hence by [D], it determines a 1-motive, which we show to be 
an extension of a direct sum of Tate structures $\Z(-1)$
by that of ($H^1$ of) an abelian variety. In the second part of the proof,
we give a direct construction of a certain 1-motive, using the Zariski
topology on $X$, and show that it is isogenous to the earlier one. The
Theorem will be an immediate corollary.

In a future work, we hope to use the second, algebraically defined
1-motive to also obtain the analogue of the Tate conjecture in our
situation, which would similarly characterize the $\Z_{\ell}$-span of the 
classes of Cartier divisors in $H^2_{\rm et\/}(X_{\bar{K}},\Z_{\ell}(1))$,
for a normal projective variety $X$ over a number field $K$.

In another direction, our result suggests a question analogous to the
Hodge conjecture. Let $X$ be a normal projective variety over $\C$, and
$\alpha:X_{an}\to X$ be the obvious continuous map, leading to a Leray
spectral sequence 
\[E_2^{p,q}=H^p(X,R^q\alpha_*\Q_{X_{an}})\implies H^{p+q}(X_{an},\Q),\]
with an induced decreasing {\em Leray filtration}
$\{L^pH^n(X_{an},\Q)\}_{p\geq 0}$ on each cohomology group
$H^n(X_{an},\Q)$. Let 
\[{\rm Hg}^p(X)=L^pH^{2p}(X_{an},\Q)\cap F^pH^{2p}(X_{an},\C).\]
Is ${\rm Hg}^p(X)$ the image of the $p$-th Chern class
map $c_p:K_0(X)\tensor\Q\to H^{2p}(X,\Q)$? Note that this does not hold
without some hypothesis like (at least) normality; for example, Bloch's
letter to Jannsen [J,Appendix A] gives a counterexample. On the other
hand, results of Collino [Co] imply it when $X$ has a unique singular
point. Note also that, unlike the standard Hodge conjecture, the positive
answer (our Theorem above) for divisors does not automatically imply a
positive answer for the case of 1-cycles, since we do not have Poincar\'e
duality. 

This work formed part of the first author's Ph.D. thesis, written at
the Tata Institute of Fundamental Research, Mumbai, submitted to the
Mumbai University in March, 1997.

\section{Some preliminaries}
\subsection{Constructible sheaves}
We will need below some technical results on constructible sheaves on
complex algebraic varieties. We begin by recalling the appropriate
definitions, the first from [V] and the second from [BS2].  

\begin{defn} Let $X$ be a complex algebraic variety. We say that a sheaf
${\cal F}$ of abelian groups on the analytic space $X_{an}$ is
(algebraically) {\em $\Z$-constructible} if there is a finite
decomposition $X=\cup_{i\in I} X_i$, where each $X_i$ is irreducible and
Zariski closed in $X$, such that if $U_i = X_i - \cup_{X_j \propsubset
X_i}X_{j}$, then each $U_i$ is non-singular, $X$ is the disjoint union of
the $U_{i}$, and ${\cal F}|_{(U_i)_{an}}$ is a locally constant sheaf
whose fibre is a finitely generated group. We call any such collection of
subsets $\{X_i\}_{i\in I}$ an {\em admissible family} of subsets for
$\cF$.  \end{defn}

\begin{defn} A sheaf ${\cal G}$ on a scheme $X$ (over an algebraically
closed field $k$, say) is said to be {\em $\Z$-constructible for the
Zariski topology} if we can express $X$ as a finite union $X=\cup X_i$,
where $X_i \subset X$ are Zariski closed, such that if $U_{i}= X_{i} -
\cup_{X_j \propsubset X_i} X_{j}$, then each $U_i$ is non-singular, $X$ is
the disjoint union of the $U_{i}$, and ${\cal G}|_{U_i}$ is a {\em
constant sheaf} associated to a finitely generated abelian group. We call
any such collection of subsets $\{X_i\}_{i\in I}$ an {\em admissible
family} for $\cG$. 
\end{defn}

\begin{rmk}
We note that in the cited works, it is not required that the ``open
strata'' $U_i$ are non-singular, but this may clearly be assumed as well
without loss of generality, by refining any given stratification which has
all the remaining properties.
\end{rmk}

Note that if $\{X_i\}_{i\in I}$ is an admissible family of subsets for a
$\Z$-constructible sheaf in either of the senses above, then there is a
natural partial order on the index set $I$ given by $j\leq i\iff X_j\subset
X_i$. Then we clearly have $U_i=X_i-\cup_{j<i}X_j$. Note that $U_i=X_i$
precisely when $i$ is a minimal element of $I$ with respect to the partial
order. 

We recall the following basic result from [V], which is made use of
below. 
\begin{thm}\label{Verdier} If $f:Y \to X$ is a
morphism of $\C$-varieties and ${\cal F}$ is a $\Z$-constructible sheaf
on $Y_{an},$ then $R^{i}f_{*}{\cal F}$ is a $\Z$-constructible sheaf on
$X_{an}$.
\end{thm}

We also need a certain general sheaf-theoretic result, which is presumably
well-known, but for which we do not know a reference.  Let $X$ be a
topological space, $\{U_i\}_{i\in I}$ any finite collection of locally
closed subsets of $X$ which stratify $X$ (\ie $X$ is the disjoint union of
the $U_i$, and for each $i$, the closure $X_i:=\bar{U_i}$ is again a union
of some $U_j$). Let $\leq $ denote the obvious partial order on $I$,
given by $i\leq j$ \iff $X_i\subset X_j$. 

Let $f_i:U_i\to X$ be the inclusion. If $\cF$ is a sheaf of abelian groups
on $X$, and $\{i_0\leq i_1\leq\cdots\leq i_p\}$ is a $p$-chain in $I$, let
\[\cF_{i_0i_1\cdots 
i_p}:=(f_{i_0})_*f_{i_0}^{-1}(f_{i_1})_*f_{i_1}^{-1}\cdots
(f_{i_p})_*f_{i_p}^{-1}\cF.\]
Note that the sheaves $\cF_{i_0\cdots i_p}$ define a cosimplicial sheaf on
the simplicial space $X\times N(I)$, where $N(I)$ denotes the nerve of $I$
(regarded as a discrete simplicial space). The augmentation $X\times
N(I)\to X$ gives rise to a complex of sheaves on $X$
\[0\to \cF\to \bigoplus_{i\in I}\cF_i\to\bigoplus_{\{i_0\leq
i_1\}\in N_1(I)}\cF_{i_0i_1}\to\cdots\to \bigoplus_{\{i_0\leq\cdots\leq
i_p\}\in N_p(I)}\cF_{i_0\cdots i_p}\to\cdots\hspace{1cm}\cdots\;\;(*)\]

\begin{lemma}\label{resol}
The above complex $(*)$ is a resolution of $\cF$.
\end{lemma}
\begin{proof}
If $x\in U_i$, then taking the stalks at $x$, we have an associated
cosimplicial abelian group $(\cF_{i_0i_1\cdots i_p})_x$, and
a  corresponding augmented complex. Clearly $(\cF_{i_0i_1\cdots i_p})_x=0$
unless $i\leq i_0$. Since the partially ordered subset $(I\geq i)=\{j\in
I\mid i\leq j\}$ has a minimal element, one sees easily that
the stalk complex at $x$ is contractible (note that if $x\in U_i$, and
$\sigma=\{i_0\leq\cdots\leq i_p\}$ is a $p$-simplex in the nerve of
$(I\geq
i$), the stalks at $x$ of $\cF_{i_0\cdots i_p}$ and $\cF_{ii_0\cdots
i_p}$ are naturally isomorphic, where $\{i\leq i_0\leq\cdots\leq i_p\}$
is the cone over $\sigma$ with vertex $i$).  
\end{proof}
\begin{rmk}
In case $\cF$ is  $\Z$-constructible for the Zariski
topology on a scheme $X$, and $\{X_i\}$ is an admissible family for $\cF$,
such that $\cF\mid_{U_i}$ is the constant sheaf associated to $A_i$, then
$\cF_{i_0i_1\cdots i_p}$ is just the constant sheaf $(A_{i_p})_{X_{i_0}}$. 
In particular, for a $\Z$-constructible sheaf in the Zariski topology, we
obtain a {\em flasque resolution}.
\end{rmk}

The key technical result of this section is the following.
\begin{lemma}\label{const1}
Let $\cA=A_{X_{an}}$ be a constant sheaf on a complex algebraic variety $X$,
and let $\cG$ be a $\Z$-constructible sheaf on $X_{an}$. Let $f:\cA\to
\cG$ be a sheaf homomorphism, and take $\cF={\rm image}\, f$. Let
$a:X_{an}\to X$ be the natural continuous map from the analytic space
$X_{an}$ to $X$, which is the identity on points. Then we have the
following. 
\begin{points}
\item $a_*\cF$ is a constructible sheaf on $X$ for the Zariski topology.
\item The natural map $a^{-1}a_*\cF\to\cF$ is an isomorphism, and the
natural map $a_*\cA\to a_*\cF$ is surjective, \ie $a_*\cF$ is a quotient
of the constant sheaf on $X$ associated to the abelian group $A$. 
\item Let $\{X_i\}_{i\in I}$ be an admissible family of subsets for $\cG$.
Then it is also admissible for $\cF$, and for $a_*\cF$. 
There is an exact sequence 
\[0\to H^0(X_{an},\cF)\to \bigoplus_{i\in
I}H^0((U_i)_{an},\cF\mid_{(U_i)_{an}})\to
\bigoplus_{\begin{array}{c}i\leq j\\
i,j\in I\end{array}}H^0((U_{j})_{an},\cF\mid_{(U_{j})_{an}})
\]  
\end{points}
\end{lemma}
\begin{proof} We first claim that if $U\subset X$ is an irreducible
(Zariski) locally closed subset such that $\cG\mid_{U_{an}}$ is locally
free, then $\cF\mid_{U_{an}}$ is a constant sheaf associated to a finitely
generated abelian group, which is a quotient of $A$. Indeed,
$\cG\mid_{U_{an}}$ corresponds to a representation of the fundamental
group of $U_{an}$ (with respect to any convenient base point), while
$A_{U_{an}}$ corresponds to the trivial representation. The sheaf map
$f\mid_{U_{an}}$ is then a morphism of local systems, whose image
$\cF\mid_{U_{an}}$ is clearly a trivial (\ie constant) local subsystem of 
$\cG\mid_{U_{an}}$.

Now let $\{X_i\}_{i\in I}$ be admissible for $\cG$. As observed above,
$\cF\mid_{(U_i)_{an}}$ is constant for each $i$, and so $\{X_i\}_{i\in I}$
is also admissible for $\cF$. From lemma~2 of [BS2], it follows that
$a_*\cF$ is $\Z$-constructible for the Zariski topology.

>From the beginning of the exact sequence $(*)$ of lemma~\ref{resol} (for
$\cF$ on $X_{an}$) we have inclusions
\[\cF\into\bigoplus_{i\in I}\cF_i,\;\;a_*\cF\into\bigoplus_{i\in
I}a_*\cF_i,\;\;a^{-1}a_*\cF\into\bigoplus_{i\in I}a^{-1}a_*\cF_i. \]
We see at once from the definitions that $a_*\cF_i$ is (the direct
image on $X$ of) a constant sheaf on $X_i$, for
each $i$, and the natural sheaf map $a^{-1}a_*\cF_i\to\cF_i$ is injective.
Since $\cA$ is a constant sheaf, we also have that $a^{-1}a_*\cA\to \cA$
is an isomorphism. Now from the commutative diagram
\[\diagram 
a^{-1}a_*\cA \dto_{\cong}\rto & a^{-1}a_*\cF \rto|<<\ahook \dto|<<\ahook
& \bigoplus_{i\in I} a^{-1}a_*\cF_i \dto|<<\ahook \\
\cA\rto|>>\tip & \cF \rto|<<\ahook & \bigoplus_{i\in I}\cF_i
\enddiagram
\]
we deduce that $a^{-1}a_*\cF\to \cF$ is an isomorphism, and that the
natural map $a^{-1}a_*\cA\to a^{-1}a_*\cF$ is surjective. This implies
that $a_*\cA\to a_*\cF$ is surjective as well, and that $\{X_i\}$ is 
admissible for $a_*\cF$. The exact sequence in (iii) of the lemma is 
obtained from the resolution of lemma~\ref{resol} for $a_*\cF$.
\end{proof}

\subsection{A homological lemma}

We prove here an abstract homological lemma (lemma~\ref{[PS]}) which is a
variant of a lemma in [PS], which we will need later. The lemma is
formulated and proved with abelian groups, but a similar argument yields
it in an arbitrary abelian category. 
Suppose we have the following $9$-diagram, in the category
of complexes of abelian groups, with exact rows and columns. 
$$\diagram & 0 \dto & 0 \dto & 0 \dto \\
0 \rto & C_{11}^{\d} \rto \dto & C_{12}^{\d} \rto 
\dto & C_{13}^{\d} \rto \dto & 0 \\
0 \rto & C_{21}^{\d} \rto \dto & C_{22}^{\d} \rto
\dto & C_{23}^{\d} \rto \dto & 0 \\    
0 \rto & C_{31}^{\d} \rto \dto & C_{32}^{\d} \rto \dto &
C_{33}^{\d} \rto \dto & 0 \\
& 0 & 0 & 0 \enddiagram$$    

Applying the cohomology functor we get an infinite double sequence with
exact rows and columns as shown below:

$$\diagram H^{i-2}(C_{22}^{\d}) \rto \dto & H^{i-2}(C_{23}^{\d})
\rto
\dto & H^{i-1}(C_{21}^{\d}) \rto \dto & H^{i-1}(C_{22}^{\d}) \rto
\dto & H^{i-1}(C_{23}^{\d}) \dto \rto &  \\ 
H^{i-2}(C_{32}^{\d}) \rto \dto & H^{i-2}(C_{33}^{\d}) \rto
\dto & H^{i-1}(C_{31}^{\d}) \rto \dto & H^{i-1}(C_{32}^{\d}) \rto 
\dto & H^{i-1}(C_{33}^{\d}) \dto \rto & \\ 
H^{i-1}(C_{12}^{\d}) \rto \dto & H^{i-1}(C_{13}^{\d}) \rto
\dto & H^{i}(C_{11}^{\d}) \rto \dto & H^{i}(C_{12}^{\d}) \rto \dto &
H^{i}(C_{13}^{\d}) \dto \rto & \\  
H^{i-1}(C_{22}^{\d}) \rto \dto & H^{i-1}(C_{23}^{\d}) \rto
\dto & H^{i}(C_{21}^{\d}) \rto \dto & H^{i}(C_{22}^{\d}) \rto \dto &
H^{i}(C_{23}^{\d}) \dto \rto & \\  
H^{i-1}(C_{32}^{\d}) \dto \rto & H^{i-1}(C_{33}^{\d})  
\dto \rto & H^{i}(C_{31}^{\d}) \dto \rto & H^{i}(C_{32}^{\d}) \dto
\rto & H^{i}(C_{33}^{\d}) \rto \dto & \\
& & & & \enddiagram$$

Suppose now that we have an element $\alpha \in H^{i}(C^{\d}_{rs})$ 
(say for example $\alpha \in H^{i}(C^{\d}_{22})$ in the above
diagram) such that $\alpha \mapsto 0$ under both the maps with domain
$H^i(C^{\d}_{rs})$. We can then do a diagram chase in the above
cohomology diagram in the following way. Suppose $\alpha\in
H^{i}(C^{\d}_{22})$; arbitrarily choose lifts
$\beta_{1} \in H^{i}(C^{\d}_{21})$ and
$\beta_{2} \in H^{i}(C^{\d}_{12})$ lifting $\alpha.$ 
Let $\beta_{1} \mapsto \gamma_{1} \in H^{i}(C^{\d}_{31})$ and let
$\beta_{2} \mapsto \gamma_{2} \in H^{i}(C^{\d}_{13}).$ Then since 
$\gamma_{1} \mapsto 0 \in H^{i}(C^{\d}_{32})$ and $\gamma_{2}
\mapsto 0 \in H^{i}(C^{\d}_{23})$, there exist $\delta_{1}$ and
$\delta_{2}$, both in $H^{i-1}(C^{\d}_{33})$, lifting
$\gamma_{1}$ and $\gamma_{2}$ respectively. 

We can do a similar diagram chase beginning with an element $\alpha\in
H^i(C^{\d}_{rs})$, for arbitrary $i,r,s$, and end up with two elements
$\delta_1,\delta_2$ in the same group $H^j(C^{\d}_{r+1\;s+1})$, where we
read the subscripts modulo 3, and $j$ is either $i-1$, $i$ or $i+1$,
depending on $(r,s)$ (we end up at  the two places  in the diagram which 
have the same entry, and are each 1 `knight's move' away from the starting point). 

Let $\bar{H^j}(C^{\d}_{r+1\;s+1})$ denote the quotient of
$H^j(C^{\d}_{r+1\;s+1})$ by the subgroup generated by the images of the
two maps in the large commutative cohomology diagram with range
$H^j(C^{\d}_{r+1\;s+1})$. For example, 
\[\bar{H^{i-1}}(C^{\d}_{3,3})=\frac{H^{i-1}(C^{\d}_{33})}{{\rm 
image}\,H^{i-1}(C^{\d}_{23})+{\rm image}\,H^{i-1}(C^{\d}_{32})}.\]

\begin{lemma} \label{[PS]} With the notation as above, we have 
$$(\delta_{1} - \delta_{2})\mapsto 0 \in \bar{H^{j}}(C^{\d}_{r+1\;s+1}).$$
\end{lemma}
\begin{proof} We first note that, by an argument with mapping cones and
cylinders (rotating the distinguished triangles in the 9-diagram), we may
assume that $\alpha \in H^{i}(C^{\d}_{22})$ without loss
of generality. For such an $\alpha$ the analogous result  for the
cohomology diagram arising from a 9-diagram in the category of sheaves
has been proved by Parimala and Srinivas [PS, Sec 3]. The proof of this
lemma is entirely analogous: regarding the given 9-diagram as a (bounded)
double complex of complexes, one considers the total complex, which is a
5-term exact sequence of complexes, say
\[0\to \cC_0\to\cC_1\to\cC_2\to\cC_3\to\cC_4\to 0.\]
Regarding this again as a double complex, there is a spectral sequence
\[E_1^{r,s}=H^s(\cC_r)\implies H^{r+s}(Tot(\cC_{\d}))\]
(the limit is in fact 0). Then the conclusion of the lemma is
interpreted as giving two (equivalent) ways of computing the differential
$E_2^{2,i}\to E_2^{4,i-1}$.
\end{proof}

\begin{rmk} An analogue of lemma~\ref{[PS]} can be formulated for a
9-diagram in the derived category of abelian groups
$$\diagram 
 C_{11} \rto \dto & C_{12} \rto 
\dto & C_{13} \rto \dto & C_{11}[1]\dto \\
 C_{21} \rto \dto & C_{22} \rto
\dto & C_{23} \rto \dto & C_{21}[1] \dto\\    
C_{31} \rto \dto & C_{32} \rto \dto &
C_{33} \rto \dto & C_{31}[1]\dto \\
C_{11}[1]\rto& C_{12}[1]\rto & C_{13}[1]\rto & C_{11}[2] \enddiagram$$    
where the rows and columns are distinguished triangles, and where the
cohomology diagram considered earlier is replaced by the diagram obtained
by applying any abelian group valued cohomological functor (of course a
still more general formulation is also possible). This is {\em false};
O.~Gabber has kindly shown us a counterexample. 
\end{rmk}

\begin{rmk} A version of the above lemma~\ref{[PS]} also appears in a
letter from U. Jannsen to B. Gross.
\end{rmk}

\section{A short exact sequence of mixed Hodge structures}
In this section we make an analysis of the mixed Hodge structure on
\[H^1(X,\ccH^1_X)=\mbox{ subgroup of Zariski locally trivial elements in
$H^2(X_{an},\Z)$}.\] 
Our goal is to describe it as an extension of a direct sum of Tate Hodge
structures $\Z(-1)$ by a polarizable pure Hodge structure of weight 1. 

Let $X$ be our given normal projective variety over $\C$. Let $Y$ be a
resolution of singularities of $X$, and let $Y_{an}$ be the 
associated analytic space of $Y$. We have the following commutative
diagram
$$\diagram Y_{an} \rto^{a^{Y}} \dto^{\pi^{an}} & Y \dto^{\pi} \\
X_{an} \rto^{a^{X}} & X \enddiagram$$

The Leray spectral sequence for the constant sheaf $\Z=\Z_{Y_{an}}$ and
the map $\pi^{an}: Y_{an} \to X_{an}$ leads to an exact sequence:
\begin{equation}\label{1} 0 \to H^1 (X_{an}, \Z) \to H^1 (Y_{an},\Z)  \to 
H^0(X_{an}, R^1 \pi^{an}_{*}\Z) \to H^2 (X_{an},\Z) \to H^{2}(Y_{an}, \Z)
\end{equation}
Note that since $X$ is normal, we have $\pi^{an}_{*}\Z \cong \Z$. 

Define a new sheaf ${\cal F}_{\Z}$ on $X_{an}$ by 
\begin{equation}\label{FZ}
{\cal F}_{\Z} = {\rm image}\,(H^1 (Y_{an}, \Z)_{X_{an}} \to
R^{1}\pi^{an}_{*}\Z).
\end{equation}
Here by $H^{1}(Y_{an},\Z)_{X_{an}}$ we mean the constant sheaf on $X_{an}$
associated to the group $H^1 (Y_{an}, \Z),$ and the map on sheaves is
induced at the level of presheaves by the restriction map on cohomology
$H^{1}(Y_{an},\Z) \to H^{1}((\pi^{an})^{-1}(U_{an}),\Z)$ where $U_{an}
\subset X_{an}$ is open. By taking global sections we have the following
commutative diagram,
$$\diagram H^{1}(Y_{an},\Z) \rto \drto & H^{0}(X_{an},{\cal F}_{\Z}) 
\dto|<<\ahook \\ & H^{0}(X_{an},R^{1}\pi^{an}_{*}\Z) \enddiagram$$

Hence we have an inclusion,
$$0 \to \frac{H^{0}(X_{an},{\cal F}_{\Z})}{{\rm Im}(H^{1}(Y_{an},\Z))}
\to \frac{H^{0}(X_{an},R^{1}\pi^{an}_{*}\Z)}{{\rm Im}(H^{1}(Y_{an},\Z))}
= Ker(H^{2}(X_{an},\Z) \to H^{2}(Y_{an},\Z))$$
where the last equality is due to the above exact sequence (\ref{1}) of 
low degree terms of the Leray spectral sequence. 

Note that $\cF_{\Z}$ satisfies the hypotheses of lemma~\ref{const1}, with
$A=H^1(Y_{an},\Z)$ and $\cG=R^1\pi^{an}_*\Z$ (the latter is algebraically
$\Z$-constructible by theorem~\ref{Verdier}). Hence the following
properties hold. 
\begin{points}
\item $\cF_{\Z}$ is algebraically $\Z$-constructible.
\item  $a^X_*\cF_{\Z}:=\cG_{\Z}$ is $\Z$-constructible for the Zariski
topology, and $(a^X)^*\cG_{\Z}\cong\cF_{\Z}$.
\item  The natural sheaf map 
\begin{equation}\label{eqsur}
H^1(Y_{an},\Z)_X\to a^X_*\cF_{\Z}
\end{equation}
is surjective. 
\end{points}

\begin{lemma}\label{incl}
 $\displaystyle{\frac{H^{0}(X_{an},{\cal
F}_{\Z})}{Im(H^{1}(Y_{an},\Z))} 
\subset H^{1}(X,{\cal H}^{1}_{X})}.$
\end{lemma}
\begin{proof} Let $\alpha \in H^{0}(X_{an},{\cal
F}_{\Z})=H^0(X,a^X_*\cF_{\Z})$. Then by (\ref{eqsur}), 
there exists a Zariski open cover
$\{U_i\}$ of $X$ such that $\alpha|_{(U_i)_{an}} = Im(\beta_{i})$ where
$\beta_{i} \in H^{1}(Y_{an},\Z).$ Therefore $\alpha|_{(U_{i})_{an}} \to 0
\in H^{2}((U_i)_{an},\Z)$ as shown in the commutative diagram below (where
$U$ stands for any of the $U_i$)
$$\diagram H^{1}(Y_{an},\Z) \dto \rto & H^{0}(X_{an},{\cal F}_{\sZ}) \rto 
\dto & H^{2}(X_{an},\Z) \dto \\
H^{1}((\pi^{an})^{-1}(U_{an}),\Z) \rto & H^{0}(U_{an},{\cal F}_{\sZ}) 
\rto & H^{2}(U_{an}\Z) \enddiagram$$
This finishes the proof of the lemma. 
\end{proof}

If $\{X_i\}_{i\in I}$ is an admissible family of subsets
for the constructible sheaf $R^1\pi^{an}_*\Z$ on $X_{an}$, then
(lemma~\ref{const1}) it is also an admissible family for $\cF_{\sZ}$ and
for $a^X_*\cF_{\sZ}=\cG_{\sZ}$. We fix such an admissible family once and
for all, and fix base points $x_i\in U_i$ with corresponding reduced
fibers $F_i=\pi^{-1}(x_i)_{red}$. Let
$F=\cup_iF_i=\pi^{-1}(\{x_i\mid i\in I\})$. 

By the proper base change theorem, the stalk
$(R^1\pi^{an}_*\Z)_{x_i}$ is naturally identified with
$H^1((F_i)_{an},\Z)$; thus $R^1\pi^{an}_*\Z\mid_{(U_i)_{an}}$ is a local
system with fiber $H^1((F_i)_{an},\Z)$.  Note that the stalk
$(\cF_{\sZ})_{x_i}$ has the resulting description
\begin{equation}\label{stalk}
(\cF_{\sZ})_{x_i}={\rm image}\,\left(H^1(Y_{an},\Z)\to
H^1((F_i)_{an},\Z)\right).
\end{equation}
By mixed Hodge theory [D], we deduce that $(\cF_{\sZ})_{x_i}$ naturally
supports a pure Hodge structure of weight 1, which is a quotient Hodge
structure of $H^1(Y_{an},\Z)$ (depending only on $i\in I$, and not on the 
chosen base point $x_i\in U_i$), as well as a Hodge sub-structure of
$H^1((F_i)_{an},\Z)$. Finally note also that $\cF_{\sZ}\mid_{(U_i)_{an}}$
is a constant sheaf whose fiber supports this pure Hodge structure (\ie is
the underlying lattice). 

\begin{lemma} $H^{0}(X_{an},{\cal F}_{\sZ})$ carries a pure Hodge 
structure of weight one, such that $H^{1}(Y_{an},\Z) \to
H^{0}(X_{an},{\cal F}_{\sZ})$ is a morphism of Hodge structures.  
\end{lemma}
\begin{proof} 
>From lemma~\ref{const1}(iii) there exists an  
exact sequence of abelian groups 
$$0 \to H^{0}(X_{an}, {\cal F}_{\sZ}) \to
\bigoplus_{i\in I}(\cF_{\sZ})_{x_i} \to 
\bigoplus_{i_0,i_1\in I,\;i_0\leq  
i_1}(\cF_{\sZ})_{x_{i_1}}$$
The natural surjective maps $(\cF_{\sZ})_{x_i}\to (\cF_{\sZ})_{x_j}$
(for $i\leq j$) are maps of pure Hodge structures of weight one, which are
quotients of $H^{1}(Y_{an},\Z)$. Hence $H^{0}(X_{an}, {\cal F}_{\sZ})$ is
identified with the kernel of a morphism of pure Hodge structures of
weight~1, and hence itself supports a pure Hodge structure of weight one.
Also it is clear from the construction that the composition 
\[H^1(Y_{an},\Z)\to H^0(X_{an},\cF_{\sZ})\into
\bigoplus_{i\in I}(\cF_{\sZ})_{x_i}\] 
is a direct sum of the natural quotient maps $H^1(Y_{an},\Z)\to
(\cF_{\sZ})_{x_i}$, and hence is a morphism of Hodge structures. Hence
$H^1(Y_{an},\Z)\to H^0(X_{an},\cF_{\sZ})$ is one as well.
\end{proof}

\begin{propose} \label{MHS} $H^{0}(X_{an},{\cal F}_{\sZ}) \to 
H^{2}(X_{an},\Z)$ is morphism of Hodge structures, i.e., the Hodge 
structures on $H^{0}(X_{an},{\cal F}_{\sZ})$ and $H^{2}(X_{an}, \Z)$ 
are compatible. 
\end{propose}
\begin{proof} Let $F_i=\pi^{-1}(x_i)$ as above, and let $F = \cup_{i\in I}
F_i$. We note that the natural map $H^{0}(X_{an}, R^{1}\pi_{*}^{an}\Z) \to
H^{1}(F_{an},\Z)$ is an injection (any section in the kernel must vanish
in all
stalks). This implies that the map $H^{2}(X_{an},\Z) \to
H^{2}(Y_{an},F_{an},\Z)$ (which is a morphism of mixed Hodge structures)
is injective in the following commutative diagram (here
$G_i=(\cF_{\sZ})_{x_i}$). 

$$\diagram \displaystyle{\frac{H^{0}(X_{an},{\cal
F}_{\sZ})}{Im(H^{1}(Y_{an},\Z)}}  
\rto|<<\ahook \dto|<<\ahook & 
\displaystyle{\frac{H^{0}(X_{an}, 
R^{1}\pi^{an}_{*}\Z)}{Im(H^{1}(Y_{an},\Z)}} 
\ddlto|<<\ahook \rto|<<\ahook & H^{2}(X_{an},\Z) \ddto|<<\ahook \\ 
\displaystyle{\frac{\oplus_{i}G_{i}}{Im(H^{1}(Y_{an},\Z)}}
\dto|<<\ahook \\ 
\displaystyle{\frac{\oplus_{i}H^{1}((F_{i})_{an}, 
\Z)}{Im(H^{1}(Y_{an},\Z))}} \rrto|<<\ahook & & 
H^{2}(Y_{an},F_{an},\Z) \enddiagram$$
We are done, because all the arrows in the above diagram are injections,
and the vertical arrows (on the left and right borders), as well as the
lower horizontal arrow, are morphisms of mixed Hodge structures.  
\end{proof}

Let $A \subset B$ be an inclusion of abelian groups. Let $A \subset A^s
\subset B$ denote the saturation of $A$ in $B,$ i.e., $A^{s}$ is the
smallest subgroup of $B$ containing $A$ such that
$\displaystyle{\frac{B}{A^{s}}}$ is torsion free. Let $H^{0}(X_{an}, {\cal
F}_{\sZ})^{s}$ be the saturation of $H^{0}(X_{an},{\cal F}_{\sZ})$ in
$H^{0}(X_{an}, R^{1}\pi^{an}_{*}\Z).$ 

\begin{lemma}  $\displaystyle{Ker(H^{1}(X,{\cal H}^{1}_{X}) \to
H^{2}(Y_{an},\Z)) = \frac{H^{0}(X_{an},{\cal
F}_{\Z})^{s}}{Im(H^{1}(Y_{an},\Z))}}.$ 
\end{lemma} 
\begin{proof} It is easy to see, from lemma~\ref{incl}, that 
\[Ker(H^{1}(X,{\cal H}^{1}_{X}) \to
H^{2}(Y_{an},\Z)) \supset \frac{H^{0}(X_{an},{\cal
F}_{\Z})^{s}}{Im(H^{1}(Y_{an},\Z))}.\]
We will prove, using lemma~\ref{[PS]}, that given any element 
$$\alpha \in ker(H^{1}(X,{\cal H}^{1}_{X}) \to H^{2}(Y_{an},\Z)),$$ 
and any preimage $\beta_{1} \in H^{0}(X_{an},R^{1}\pi^{an}_{*}\Z),$
some non-zero (integral) multiple of $\beta_{1}$ lies in
$H^{0}(X_{an},{\cal F}_{\Z}) \subset H^{0}(X_{an},R^{1}\pi^{an}_{*}\Z).$
This will prove the assertion of the lemma.

Since $\alpha \in H^{1}(X,{\cal H}^{1}_{X}),$ there exists a finite
Zariski open cover $\{U_{i}\}$ of $X$ such that $\alpha \mapsto 0$ in
$H^{2}((U_{i})_{an},\Z)$ for all $i.$ Let $U$ denote any one of these
$U_{i}$'s and consider again the above commutative diagram with exact
rows and columns.  

$$\hspace{-1cm}\diagram & & & H^{1}(Y_{an},\Z) \dto \\
& & H^{1}(U_{an},\Z) \rto \dto & H^{1}(\pi^{-1}(U_{an},\Z) \dto \\
& & H^{2}(X_{an},U_{an},\Z) \rto \dto & 
H^{2}(Y_{an},\pi^{-1}(U_{an}),\Z) \dto \\
H^{1}(Y_{an},\Z) \rto \dto & 
H^{0}(X_{an},R^{1}\pi^{an}_{*}\Z) \rto \dto & H^{2}(X_{an},\Z) \rto \dto 
& H^{2}(Y_{an},\Z) \dto \\
H^{1}(\pi^{-1}(U_{an}),\Z) \rto &  
\Gamma(U_{an},R^{1}\pi_{*}^{an}\Z) \rto & H^{2}(U_{an},\Z) \rto & 
H^{2}(\pi^{-1}(U_{an},\Z) \enddiagram$$

We wish to apply lemma~\ref{[PS]} to this diagram; for this, we need to
know that this diagram arises by applying the cohomology functor to a
suitable $9$-diagram in the category of complexes of abelian groups. It is
clear that the above diagram arises by applying the cohomology functor to
the following $9$-diagram, where all the objects are in the (bounded
below) derived category of sheaves of abelian groups on $X$, and the rows
and columns are exact triangles; here $K_i$ are suitable cones. 
$$\diagram K_{1} \rto \dto & K_{2} \rto \dto & K_{3} \dto \rto & \\
\Z_{X} \rto \dto & R\pi_{*}\Z_{Y} \rto \dto & C_{1} \rto \dto & \\
Rj_{*}\Z_{U} \rto \dto & Rj_{*}R\pi_{*}\Z_{\pi^{-1}(U)} \rto \dto & 
Rj_{*}C_{2} \rto \dto & \\
& & & \enddiagram$$

Applying the functor $R\Gamma(X,-)$ yields a $9$-diagram in the derived
category of abelian groups. Using Cartan-Eilenberg resolutions, this
$9$-diagram in the derived category is seen to be the image of a
$9$-diagram where all the objects are complexes of abelian groups and the
rows and columns are short exact sequences of complexes. Since the
arguments are standard, we omit the details.

Returning to our cohomology diagram, note that the relative cohomology
sequences  
$$\to H^{1}(U_{an},\Z) \to H^{2}(X_{an},U_{an},\Z) \to H^{2}(X_{an},\Z)
\to H^{2}(U_{an},\Z) \to$$
and 
$$\hspace{-1cm}\to H^{1}((\pi^{an})^{-1}(U_{an}),\Z) \to 
H^{2}(Y_{an},(\pi^{an})^{-1}(U_{an}),\Z) \to
H^{2}(Y_{an},\Z) \to H^{2}((\pi^{an})^{-1}(U_{an},\Z)$$ 
are sequences in the category of mixed Hodge structures by 
[D]~(8.3.9). 

Since $\alpha \to 0 \in H^{2}(Y_{an},\Z)$ therefore $\alpha_{\Q} \in
W_{1}H^{2}(X_{an},\Q)$ by [D],~Proposition~8.2.5. This implies, by 
[D],~Theorem~2.3.5 (\ie strictness of morphisms of mixed Hodge structures
with respect to $W$) that we can choose $\beta_{2} \in H^{2}(X_{an}, 
U_{an},\Z)$ such that 
\[(\beta_{2})_{\Q} \in W_{1}H^{2}(X_{an},U_{an},\Q),\;\;\beta_{2} \mapsto
n \alpha,\;\;n\in\Z_{>0}.\]
Let 
\[\beta_2\mapsto\gamma_{2}\in
H^{2}(Y_{an},(\pi^{an})^{-1}(U_{an}),\Z) \cong \Z(-1)^k,\]
for some $k\geq 0$, where the last isomorphism is because  $Y$ is
non-singular; then
\[(\gamma_{2})_{\Q} \in W_{1}H^{2}(Y_{an},(\pi^{an})^{-1}(U_{an}),\Q)=0,\]
\ie $\gamma_2=0$. So we can choose a preimage $\delta_{2} \in
H^{1}((\pi^{an})^{-1}(U_{an}),\Z)$ of $\gamma_{2}$ to be zero.
On the other hand, chasing the diagram the other way, we get $n\beta_{1}
\in H^{0}(X_{an},R^{1}\pi^{an}_{*}\Z)$ which lifts $n 
\alpha$, and $n\beta\mapsto n\gamma_{1} \in 
H^{0}(U_{an},R^{1}\pi^{an}_{*}\Z)$; now take a lift $n\delta_{1} \in
H^{1}((\pi^{an})^{-1}(U_{an},\Z).$ of $n\gamma_1$.

By lemma~\ref{[PS]}, we know that $n\delta_{1} \equiv \delta_{2}=0$
modulo the images of $H^{1}(Y_{an},\Z)$ and $H^{1}(U_{an},\Z).$ Therefore 
\[n\gamma_{1} \in Im(H^{1}(Y_{an},\Z)\to
H^{0}(U_{an},R^{1}\pi^{an}_{*}\Z)).\]
This proves that  
$n\beta_{1}|_{U_{an}}$ comes from $H^{1}(Y_{an},\Z).$ 
Since $X$ has a finite cover by such open sets $U,$ we see that $\beta_{1}
\in H^{0}(X_{an},{\cal F}_{\sZ})^{s}.$ 
\end{proof}

\begin{cor}\label{corexseq}
There exists a short exact sequence of mixed Hodge structures
$$0 \to
\frac{H^{0}(X_{an},{\cal F}_{\sZ})^{s}}{Im(H^{1}(Y_{an},\Z))} \to
H^{1}(X,{\cal H}^{1}_{X}) \to 
Im(H^{1}(X,{\cal H}^{1}_{X}) \to H^{2}(Y_{an},\Z)) \to 0.$$ 
\end{cor}

Let $Im(H^{1}(Y_{an},\Z))^{s}$ denote the saturation of $H^{1}(Y_{an},\Z)$
in $H^{0}(X_{an},{\cal F}_{\sZ})^{s}.$ Then
\[\frac{Im(H^{1}(Y_{an},\Z))^{s}}{Im(H^{1}(Y_{an},\Z)} =
\left(\frac{H^{0}(X_{an},{\cal F}_{\sZ})^{s}}
{Im(H^{1}(Y_{an},\Z)}\right)_{\rm torsion}.\] Since 
$NS(X) \subset F^{1}H^{2}(X_{an},\Z)$ it follows that it has finite
intersection with $\displaystyle{\frac{H^{0}(X_{an},{\cal
F}_{\Z})^{s}}{Im(H^{1}(Y_{an},\Z))}}$ which is a pure Hodge structure of
weight one. On the other hand, $H^{2}(X_{an},\Z)_{\rm torsion} \subset
NS(X)$ from the exponential sequence. Thus, 
$$\left(
\frac{H^{0}(X_{an},{\cal F}_{\sZ})^{s}}{Im(H^{1}(Y_{an},\Z)}
\right)_{\rm torsion} = NS(X) \cap
\frac{H^{0}(X_{an},{\cal F}_{\sZ})^{s}}{Im(H^{1}(Y_{an},\Z)}.$$
Hence we get the exact sequence
$$0 \to
\displaystyle{\frac{H^{0}(X_{an},{\cal
F}_{\sZ})^{s}}{Im(H^{1}(Y_{an},\Z))^{s}}} \to
\displaystyle{\frac{H^{1}(X,{\cal H}^{1}_{X})}{NS(X)}}
\stackrel{f}{\to} \displaystyle{\frac{Im(H^{1}(X,{\cal H}^{1}_{X}) \to
H^{2}(Y_{an},\Z))}{Im(NS(X))}} \to 0,\hspace{2mm}(+)$$ 
It is clear that the third term is pure of type $(1,1)$ as it lies
inside $\displaystyle{\frac{NS(Y)}{Im(NS(X))}}.$ 

Let $A =  \displaystyle{\frac{Im(H^{1}(X,{\cal H}^{1}_{X}) \to 
H^{2}(Y_{an},\Z))}{Im(NS(X))}}$ and let 
$$M = f^{-1}(A_{\rm torsion}).$$ 
We then have a short exact sequence of mixed Hodge structures
$$0 \to M \to \frac{H^{1}(X,{\cal H}^{1}_{X})}{NS(X)} \to
\frac{A}{A_{\rm torsion}} \cong \Z(-1)^{r} \to 0,\hspace{2mm}(++)$$
where the third term is free of rank $r$ and pure of type $(1,1),$ and $M$
is a pure Hodge structure of weight $1.$ Further, all of the underlying
abelian groups are free. 

We recall some facts about extensions of mixed Hodge structures (see [C],
for example). Let $H$ be a finitely generated abelian 
group which supports a pure  Hodge structure of weight one, and $G$ a
finitely generated abelian group, regarded as a pure Hodge structure of
type $(0,0)$. Then there is a natural identification of the abelian group
$\Ext^{1}_{\bf MHS}(G(-1),H)$ of extensions in the category ${\bf MHS}$ of
mixed Hodge structures with the abelian group $\Hom(G,J(H))$, where
$$J(H) = J^{1}(H) = \frac{H_{\sC}}{F^{1}H_{\sC}{}+{}Im(H)};$$ 
here $H_{\sC} = H \tensor_{\sZ}\C$ and $F$ gives the Hodge filtration.
In particular we have
\[\Ext^1_{\bf MHS}(\Z(-1),H)=J(H).\]
If 
\[0\to H\to E\to G(-1)\to 0\]
is an extension of mixed Hodge structures, let $\psi_E:G\to J(H)$ be the
corresponding homomorphism (which we call the {\em extension class map} of
$E$). This may be described as follows: there is an identification
\[\alpha:\frac{H_{\sC}}{F^1H_{\sC}}\by{\cong}\frac{E_{\sC}}{F^1E_{\sC}},\]
giving
\[\beta:J(H)=\frac{H_{\sC}}{F^1H_{\sC}+H}\by{\cong}\frac{E_{\sC}}{F^1E_{\sC}+H},\]
and $\psi_E$ is the composition
\[G\cong \frac{E}{H}\to \frac{E_{\sC}}{F^1E_{\sC}+H}
\by{\beta^{-1}}J(H).\] 
In case $G=\Z^{\oplus r}$ is free abelian, we have that $\Hom_{\bf
MHS}(\Z(-1),G)$ is naturally identified with $\ker \psi_E$. Also, if $G$
is free abelian and $H$ is polarizable, then $J(H)$ is an abelian
variety, and for an extension $E$, the homomorphism $\psi_E:G\to J(H)$ is
the 1-motive over $\C$ associated to the mixed Hodge structure $E$ by the
procedure in [D], (10.1.3). 

In particular, the sequence of mixed Hodge structures $(+)$ is an 
extension of a pure Hodge structure of type $(1,1)$ by a pure weight one
Hodge structure and hence gives rise to an extension class homomorphism, 
$$\psi: \frac{Im(H^{1}(X,{\cal H}^{1}_{X}) \to
H^{2}(Y_{an},\Z))}{NS(X)}
\to J\left(\frac{H^{0}(X_{an},{\cal 
F}_{\sZ})^{s}}{(Im(H^{1}(Y_{an},\Z)))^{s}}\right).$$ 
Similarly, the sequence of mixed Hodge structures $(++)$ gives rise to a
related homomorphism
$$\psi_{1}: \Z^{\oplus r} \cong \frac{A}{A_{tors}} \to J(M),$$ 
which is in fact a 1-motive. Note that $J(M)$ is isogenous to
$\displaystyle{J\left(\frac{H^{0}(X_{an},{\cal
F}_{\Z})^{s}}{(Im(H^{1}(Y_{an},\Z)))^{s}}\right)}$ which in turn is
isogenous to $\displaystyle{J\left(\frac{H^{0}(X_{an},{\cal
F}_{\Z})^{s}}{Im(H^{1}(Y_{an},\Z))}\right)}.$

By the above remarks, our main result Theorem~\ref{(1,1)} is
equivalent to proving {\em $\psi_{1}$ is injective}.

\section{Proof of the Main Theorem}
\subsection{Construction of a 1-motive}
The aim of this section is to directly construct a certain 1-motive over
$\C$. The proof of the Main Theorem  will be by showing that it is
isogenous to that associated to $(H^1(X,\ccH^1_X)/NS(X))\tensor\Z(1)$. 

Let $\pi : Y \to X$ be a desingularization of $X$ as before and let $U
\subset X$ be a Zariski open subset. We have an exact sequence of groups
$$Pic^{0}(Y) \to Pic(\pi^{-1}(U)) \to H^{1}(\pi^{-1}(U),{\cal H}^{1}_{Y}) 
\to 0.$$ 
We sheafify this on $X=X_{Zar}$ to get an exact sequence of
sheaves 
$$Pic^{0}(Y)_{X} \to R^{1}\pi_{*}{\cal O}^{*}_{Y} \to R^{1}\pi_{*}{\cal
H}^{1}_{Y} \to 0$$ 
where $Pic^{0}(Y)_{X}$ is the constant sheaf on $X$ associated to the
group $Pic^{0}(Y).$ Define ${\cal F}$ to be the Zariski sheaf 
$${\cal F}:=Im(Pic^{0}(Y)_{X} \to R^{1}\pi_{*}{\cal O}^{*}_{Y})$$ on $X.$
Hence we have short exact sequence 
of sheaves 
\begin{equation}\label{eq}
 0 \to {\cal F} \to R^{1}\pi_{*}{\cal O}^{*}_{Y} \to
R^{1}\pi_{*}{\cal H}^{1}_{Y} \to 0
\end{equation}

\begin{lemma} \label{compare}(1) There is an injective map
\[\mu:J(H^0(X_{an},\cF_{\sZ})^s)\to H^0(X,\cF),\]
whose image $H^{0}(X,{\cal F})^{0}$ is a subgroup of finite index, such
that the natural map $Pic^{0}(Y) \to  H^{0}(X,{\cal F})$ factors
through $\mu$. The induced map $Pic^{0}(Y) \to  H^{0}(X,{\cal F})^{0}$ is
that determined by the map on Hodge structures $H^1(Y_{an},\Z)\to
H^0(X_{an},\cF_{\sZ})^s$. Thus $H^{0}(X,{\cal F})^{0}$ is the group of
$\C$-points of an abelian variety, such that $Pic^{0}(Y) \to
H^{0}(X,{\cal F})^{0}$ is a homomorphism of abelian varieties. 

(2)\quad $\displaystyle{J\left(\frac{H^{0}(X_{an},{\cal 
F}_{\sZ})^{s}}{Im(H^{1}(Y_{an},\Z)}\right)}$ is isogenous to
$\displaystyle{\frac{H^{0}(X,{\cal F})^{0}}{Im(Pic^{0}(Y))}},$ and hence
the latter is isogenous to $J(M).$  
\end{lemma}
\begin{proof} 
Let $x \in X$ be any point and $F_{x} = \pi^{-1}(x)_{red}.$ Let $Y_x =
Spec({\cal O}_{X,x}) \times_{X} Y$ and $\displaystyle{F_{x}^{n} =
Spec\left(\frac{{\cal O}_{X,x}}{{\cal M}_{x}^{n}}\right) \times_{X}Y}$,
where ${\cal M}\subset {\cal O}_{X,x}$ is the maximal ideal.
For each $n$ we have the restriction maps $h_{n}: Pic(F^{n}_{x}) \to
Pic(F_{x}).$

We claim that the kernel of $h_{n},$ for each $n,$ is a $\C$-vector
space. To see this consider the short exact sequence of Zariski sheaves 
$$\diagram 0 \rto & {\cal I}_{F^{n}_{x}/F_{x}} \rto^{exp} & {\cal
O}^{*}_{F^{n}_{x}} \rto^{h_n} & {\cal O}^{*}_{F_{x}} \rto & 0
\enddiagram$$ 
where $exp$ denotes the exponential map, which makes sense as ${\cal
I}_{F^{n}_{x}/F_{x}}$ is nilpotent.   
Considering the associated cohomology sequence we get 
$$0 \to H^{1}(F_{x}, {\cal I}_{F^{n}_{x}|F_{x}}) \to
Pic(F^{n}_{x}) \longby{h_n} Pic(F_{x})$$
which proves the kernel of $h_{n}$, for each $n$, is a $\C$-vector space.

For each $n$ and for each $x \in X,$ we have a commutative diagram
$$\diagram Pic^{0}(Y) \rrtou ^{f_{n}} \rto \drto^{g} & Pic(Y_{x}) \dto
\rto & Pic(F^{n}_{x}) \dlto_{h_{n}} \\
& Pic(F_{x}) \enddiagram$$

Let $f_{n}$ be the composition $Pic^{0}(Y) \to Pic(Y_{x}) \to
Pic(F^{n}_{x}).$ Then, $Ker(f_{n})$ and $Ker(g)$ are both closed  
subgroups of $Pic^{0}(Y)$ hence are compact (topological) groups.
Since $Ker(h_{n})$ is a $\C$-vector space it follows that
$Ker(f_{n}) = Ker(g),$ as any continuous homomorphism from a 
compact group to a $\C$-vector space is zero (note that $Pic^{0}(Y),$
$Pic^{0}(F_{x}^{n}),$ and $Pic^{0}(F_{x})$ are isomorphic to the
corresponding analytic groups, by GAGA, and hence from the exponential
sequence carry natural topologies, such that the restriction homomorphisms
are continuous). 
  
Passing to the inverse limit we have a commutative diagram 
$$\diagram Pic^{0}(Y) \rto^{f} \drto^{g} & Pic(Y_{x}) \dto \rto|<<\ahook & 
Pic(\hat{Y_{x}}) \dlto_{h} \\ 
& Pic(F_{x}) \enddiagram$$
where $\hat{Y_{x}}$ stands for the completion of $Y_{x}$ along $F_x$.
By Grothendieck's Formal Function Theorem [H, Ch.III, Th.11.1] and the
fact that $Pic(\hat{Y_{x}}) \to \liminv{n}(Pic(F^{n}_{x}))$ is an
isomorphism [H, Ch.II, Ex.9.6], we have that $Pic(Y_{x}) \to
Pic(\hat{Y_{x}})$ is an injection. Thus it follows that $Ker(f) =
Ker(g).$ We have from the definition of ${\cal F}$ that the stalk of
${\cal F}$ at $x$, ${\cal F}_{x} = Im(Pic^{0}(Y) \to Pic(Y_{x})).$ By
our analysis so far we have proved that the natural map ${\cal F}_{x} \to
Pic(F_{x})$ is an inclusion, and it clearly factors through the the
subgroup $Pic^{0}(F_{x}).$ 

By the results of Du Bois [DB] there exists a commutative triangle  
$$\diagram H^{1}(F_{x},\C) \rto \drto & H^{1}(F_{x},{\cal O}_{F_{x}}) 
\dto^{\alpha} \\
& \displaystyle{\frac{H^{1}(F_{x},\C)}{F^{1}H^{1}(F_{x},\C)}}
\enddiagram$$
Note that $Ker(\alpha)$ is a $\C$-vector space.
Now $\alpha$ induces a map 
$$\beta : Pic^{0}(F_{x}) \cong \frac{H^{1}(F_{x},{\cal
O})}{H^{1}(F_{x},\Z)} \to J(H^{1}(F_{x},\Z)).$$ 
Thus we have a diagram
$$\diagram H^{1}(F_{x},{\cal O}_{F_{x}}) \rto^{\alpha} \dto &
\displaystyle{\frac{H^{1}(F_{x},\sC)}{F^{1}H^{1}(F_{x},\sC)}} \dto \\
Pic^{0}(F_{x}) \rto^{\beta} & J(H^{1}(F_{x},\Z)) \enddiagram$$ 
Since $H^{1}(F_{x},\Z)$ is a mixed Hodge structure
with weights $0$ and $1$ (by [D2], as $F_{x}$ is a projective variety),
$H^{1}(F_{x},\Z)$ injects into    
$\displaystyle{\frac{H^{1}(F_{x},\C)}{F^{1}H^{1}(F_{x},\C)}}$.
Thus it is clear from the above diagram $Ker(\beta) = Ker(\alpha)$ and
so $Ker(\beta)$ is also a $\C$-vector
space. Hence the composite ${\cal F}_{x} \to Pic^{0}(F_{x}) \to
J(H^{1}(F_{x},\Z))$ is injective, as ${\cal F}_{x}$ is a compact group, 
from its definition. 

Let ${\cal F}_{\sZ,x} = Im(H^{1}(Y,\Z) \to H^{1}(F_{x},\Z))$ be the stalk 
of ${\cal F}_{\sZ}$ at $x.$ Let ${\cal F}_{\sZ,x}^{s}$ be the saturation
of ${\cal F}_{\sZ,x}$ in $H^{1}(F_{x},\Z).$ The inclusion  ${\cal
F}_{\sZ,x} \to H^{1}(F_{x},\Z)$ induces a natural map with finite kernel    
$J({\cal F}_{\sZ,x}) \to J(H^{1}(F_{x},\Z)).$ In fact this map factors
as  
$$\diagram J({\cal F}_{\sZ,x}) \rto|>>\tip & J({\cal F}^{s}_{\sZ,x})   
\rto & J(H^{1}(F_{x},\Z)). \enddiagram$$
The second map is an inclusion and $J({\cal F}^{s}_{\sZ,x})$ is the
image of $J({\cal F}_{\sZ,x})$ in $J(H^{1}(F_{x},\Z)).$ 

We thus have a commutative diagram with surjective and injective maps
as follows (where we identify $Pic^{0}(Y)$ with $J(H^{1}(Y_{an},\Z)).$    

$$\diagram Pic^{0}(Y) \rto|>>\tip \dto|>>\tip & J({\cal F}_{\sZ,x}) 
\dto|>>\tip \rto & J({\cal F}_{\sZ,x}^{s}) \dlto|<<\ahook \\
{\cal F}_{x} \rto|<<\ahook & J(H^{1}(F_{x},\Z)) \enddiagram$$
 
Since it is clear from the diagram that $J({\cal F}^{s}_{\sZ,x})$ and 
${\cal F}_{x}$ are both the image of $Pic^{0}(Y)$ in
$J(H^{1}(F_{x},\Z))$ it follows that $J({\cal F}^{s}_{\sZ,x}) \cong {\cal
F}_{x}.$ Therefore there exists a map $J({\cal F}_{\sZ,x}) \to {\cal
F}_{x}$ which is an isogeny. 

We had proved that the sheaf ${\cal F}_{\sZ}$ was constructible, i.e.,
constant with groups $G_{i}$ over locally closed sets $(U_{i})_{an},$ and
this data gives rise to a flasque resolution of $a_*\cF_{\sZ}$ in the
Zariski site (by lemma~\ref{const1} and lemma~\ref{resol}). It is then
clear that analogous results hold also for the sheaf ${\cal F}_{\sZ}^{s}$
where ${\cal F}_{\sZ}^{s}$ denotes the saturation of the sheaf ${\cal
F}_{\sZ}$ in $R^{1}\pi^{an}_{*}\Z$ (which is a torsion-free sheaf). Thus
the abelian varieties $J({\cal F}^{s}_{\sZ,x})$ are constant quotients of
$Pic^{0}(Y)$ over the strata $U_{i},$ hence so are ${\cal F}_{x}.$ This
proves that the sheaf ${\cal F}$ is a constructible sheaf on $X$ for the
Zariski topology, with admissible family $\{X_{i}\},$ and further (by
lemma~\ref{resol}) ${\cal F}$ has a flasque resolution similar to
$a_*{\cal F}_{\sZ}.$

Taking global sections of the flasque resolution of $a_*{\cal
F}^{s}_{\sZ}$, we get an exact sequence 
$$0 \to H^{0}(X_{an},{\cal F}^{s}_{\sZ}) \to \bigoplus_{i}{\cal
F}^{s}_{\sZ,x_{i}} \to \bigoplus_{i<j}{\cal F}^{s}_{\sZ,x_{j}}.$$  
Also it is clear from the definitions that
$H^{0}(X_{an},{\cal F}^{s}_{\sZ})=H^{0}(X_{an},{\cal F}_{\sZ})^s$.
Applying $J$ on all the terms, we obtain a complex 
$$0 \to J(H^{0}(X_{an},{\cal F}^{s}_{\sZ})) \to \oplus_{i}J({\cal
F}^{s}_{\sZ,x_{i}}) \to \oplus_{i<j}J({\cal F}^{s}_{\sZ,x_{j}})$$  
This complex is exact on the left and has   
finite homology in the middle, since $H^{0}(\oplus_{i<j}{\cal
F}^{s}_{\sZ,x_{j}})$ is torsion-free. 
Similarly taking global sections of the flasque resolution of ${\cal
F}$ we get an exact sequence
$$0 \to H^{0}(X,{\cal F}) \to \oplus_{i}{\cal F}_{x_{i}} \to
\oplus_{i<j}{\cal F}_{x_{j}}.$$ 

There exists a commutative diagram 
$$\diagram & J(H^{0}(X,{\cal F}_{\sZ})^{s}) \ddotted_{\mu} \rto &
\oplus_{i}J({\cal 
F}_{\sZ,x_{i}}^{s}) \rto \ddouble & \oplus_{i<j}J({\cal 
F}_{\sZ,x_{j}}^{s}) \ddouble \\ 
0 \rto & H^{0}(X,{\cal F}) \rto & \oplus_{i}{\cal F}_{x_{i}} \rto & 
\oplus_{i<j}{\cal F}_{x_j} \enddiagram$$
where the two vertical arrows are isomorphisms. Hence the dotted arrow 
$\mu$ exists, and is an inclusion with finite cokernel. 

Define 
$$H^{0}(X,{\cal F})^{0} = Im(J(H^{0}(X,{\cal F}_{\sZ})^{s}).$$
Clearly this is an abelian variety, and there is an isogeny 
$J(H^{0}(X,{\cal F}_{\sZ})) \to H^{0}(X,{\cal F})^{0}.$ Also, by
construction, the natural
map $\Pic^0(Y)\to H^0(X,\cF)$ clearly factors through the map
\[\Pic^0(Y)=J(H^1(Y_{an},\Z))\to J(H^{0}(X,{\cal F}_{\sZ})).\]

Thus we have 
an isogeny $\displaystyle{\frac{J(H^{0}(X_{an},{\cal  
F}_{\sZ}))}{Im(Pic^{0}(Y))} \to \frac{H^{0}(X, {\cal
F})^{0}}{Im(Pic^{0}(Y))}}.$                                         
We finally note that there exists an isogeny 
$$\frac{J(H^{0}(X_{an},{\cal F}_{\sZ}))}{Im(Pic^{0}(Y)} \to
J\left(\frac{H^{0}(X_{an},{\cal F}_{\sZ})}{Im(H^{1}(Y_{an},\Z)}\right)$$ 
since $J(H^{1}(Y_{an},\Z)) \cong Pic^{0}(Y).$ 
                                                                       
This finishes the proof that $\displaystyle{\frac{H^{0}(X,{\cal
F})^{0}}{Im(Pic^{0}(Y))}}$ and $\displaystyle{J\left  
(\frac{H^{0}(X_{an},{\cal
F}_{\sZ})}{Im(H^{1}(Y_{an},\Z)}\right)}$ are isogenous. 
\end{proof}

We can now construct a 1-motive, as follows. Since $X$ is normal, we have
that $\pi_{*}{\cal O}_{Y}= {\cal O}_{X}$, and so we have an exact sequence 
$$0 \to Pic(X) \to Pic(Y) \to H^{0}(X,R^{1}\pi_{*}{\cal O}^{*}_{Y}).$$  
This induces another exact sequence 
\[NS(X) \to NS(Y) \to \frac{H^{0}(X,R^{1}\pi_{*}{\cal
O}^{*}_{Y})}{Im(Pic^{0}(Y))}.\]
We thus have an {\em injective} map 
\begin{equation}\label{eqq}
\frac{NS(Y)}{Im(NS(X))} \to \frac{H^{0}(X,R^{1}\pi_{*}{\cal
O}^{*}_{Y})}{Im(Pic^{0}(Y))}\end{equation} 

\begin{lemma}\label{lemphi}
 The map (\ref{eqq}) induces an (injective) map
\[\phi:\frac{Im(H^{1}(X,{\cal H}^{1}_{X}) \to H^{2}(Y_{an},\Z))}{Im(NS(X))}
\to \frac{\Gamma(X,{\cal F})}{Im(Pic^{0}(Y))}.\]
\end{lemma}
\begin{proof}
Using the short exact sequence of sheaves (\ref{eq}), we get the following
commutative diagram, whose right column is exact,  
\[\diagram \displaystyle{\frac{Im(H^{1}(X,{\cal H}^{1}_{X}) \to
H^{2}(Y_{an},\Z))}{Im(NS(X))}} \dto|<<\ahook \rdotted^>>>{\phi} &
\displaystyle{\frac{\Gamma(X,{\cal F})}{Im(Pic^{0}(Y))}} \dto \\
\displaystyle{\frac{NS(Y)}{Im(NS(X))}} \rto|<<\ahook &
\displaystyle{\frac{\Gamma(X,R^{1}\pi_{*}{\cal O}^{*})}{Im(Pic^{0}(Y))}}
\dto \\ & \Gamma(X,R^{1}\pi_{*}{\cal H}^{1}_{Y}) \enddiagram\]
Here, we claim the dotted arrow $\phi$ exists (and is also injective)
because the composition
\[\frac{Im(H^{1}(X,{\cal H}^{1}_{X}) \to H^{2}(Y_{an},\Z))}{NS(X)} \to
H^{0}(X,R^{1}\pi_{*}{\cal H}^{1}_{Y})\]
is zero. This is obvious as this map can be described in the following
way: given the image in $H^1(Y,\ccH^1_Y)=NS(Y)$  of a Zariski locally
trivial cohomology class $\eta\in H^1(X,\ccH^1_X)$, consider a line bundle 
$L_{\eta}$ on $Y$ which represents it, then consider the line bundle
restricted to open sets $\pi^{-1}(U) \subset Y,$ 
$L_{\eta}|_{\pi^{-1}(U)},$ (where $U \subset X$ open) and take the Chern
classes of these restrictions.  These give a global section of
$R^{1}\pi_{*}{\cal H}^{1}_{Y}$ which is zero as the line bundle came from
a locally trivial cohomology class on $X.$  
\end{proof}

Since $\displaystyle{\frac{\Gamma(X,{\cal F})}{Im(Pic^{0}(Y))}}$ has a
subgroup of finite index which is an abelian variety, $\phi$ determines a
1-motive in an obvious way,
\[B\to  \frac{\Gamma(X,{\cal F})^0}{Im(Pic^{0}(Y))},\]
where $B$ is the inverse image under $\phi$ of the abelian variety. 

\begin{rmk} We do not know if $\Gamma(X,\cF)$ is itself an abelian
variety, \ie if $\Gamma(X,\cF)^0=\Gamma(X,\cF)$. 
\end{rmk}

\subsection{Comparison of the two 1-motives}

We now finish the proof of the theorem, by comparing the 1-motive
constructed above using $\phi$ with that constructed earlier, using the
extension class map $\psi$ for the mixed Hodge structure on
$H^1(X,\cHH^1_X)$.  

Recall that $\{X_i\}_{i\in I}$ is the chosen admissible family of subsets
for $R^1\pi^{an}_*\Z$, and hence for $\cF_{\sZ}$ and $\cF$ as well; recall
also the corresponding (irreducible, non-singular) locally closed strata
$\{U_i\}_{i\in I}$. Also recall the choice of points $x_i\in U_i$, and
$F_i=\pi^{-1}(x_i)$, $F=\cup_{i\in I}F_i$.

Suitably blow up $Y$ to get $f:\tilde{Y} \to Y$, with $\tilde{Y}$
non-singular projective, so that the reduced strict transform of $F$ is a
{\em smooth} possibly disconnected subvariety $\tilde{F}.$ Note that there
exists the following diagram.
$$\diagram \tilde{F} \rto|<<\ahook \dto & \tilde{Y} \dto^{f} \\
F \rto|<<\ahook & Y \enddiagram$$

Now consider the following commutative diagram which is a diagram in the
category of mixed Hodge structures 
$$\hspace{-2cm}\diagram 0 \rto &
\displaystyle{\frac{H^{0}(X_{an},{\cal
F}_{\sZ})^{s}}{Im(H^{1}(Y_{an},\Z))}}\rto \dto & H^{1}(X,{\cal
H}^{1}_{X}) \rto \dto & Im(H^{1}(X,{\cal H}^{1}_{X}) \to
H^{2}(Y_{an},\Z)) \rto \dto & 0 \\ 0 \rto &
\displaystyle{\frac{H^{1}(F_{an},\Z)}{Im(H^{1}(Y_{an},\Z))}} \rto \dto &
H^{2}(Y_{an},F_{an},\Z) \rto \dto & Im(H^{2}(Y_{an},F_{an},\Z) \to
H^{2}(Y_{an},\Z)) \dto \rto & 0 \\
0 \rto &
\displaystyle{\frac{H^{1}(\tilde{F_{an}},\Z)}{Im(H^{1}(\tilde{Y_{an}}  
,\Z))}} \rto & H^{2}(\tilde{Y_{an}},\tilde{F_{an}},\Z) \rto & 
Im(H^{2}(\tilde{Y_{an}},\tilde{F_{an}},\Z) \to H^{2}(\tilde{Y_{an}},\Z))
\rto & 0 \enddiagram$$

Let $Im(H^{1}(\tilde{Y_{an}},\Z))^{s}$ be the saturation of
$Im(H^{1}(\tilde{Y_{an}},\Z))$ in $H^{1}(\tilde{F_{an}},\Z).$ Then,
arguing as before, 
$$Im(NS(X)) \cap  
\frac{H^{1}(\tilde{F_{an}},\Z)}{Im(H^{1}(\tilde{Y_{an}},\Z))} =
\frac{Im(H^{1}(\tilde{Y_{an}},\Z))^{s}}{Im(H^{1}(\tilde{Y_{an}},\Z))}.$$
So we have a short exact sequence of mixed Hodge structures 
$$\diagram 0 \rto &
\displaystyle{\frac{H^{1}(\tilde{F_{an}},\Z)}{Im(H^{1}(\tilde{Y_{an}}
,\Z))^{s}}} \rto & 
\displaystyle{\frac{H^{2}(\tilde{Y_{an}},\tilde{F_{an}},\Z)}{Im(NS(X))}} 
\rto & 
\displaystyle{\frac{Im(H^{2}(\tilde{Y_{an}},\tilde{F_{an}},\Z) 
\to H^{2}(\tilde{Y_{an}},\Z))}{Im(NS(X))}} \rto & 0 \enddiagram$$ 
and a commutative diagram 
$$\diagram 0 \rto &
\displaystyle{\frac{H^{0}(X_{an},{\cal F}_{\sZ})^{s}}
{Im(H^{1}(Y_{an},\Z))^{s}}} \rto \dto &
\displaystyle{\frac{H^{1}(X,{\cal H}_{X}^{1})}{NS(X)}} \rto \dto &
\displaystyle{\frac{Im(H^{1}(X,{\cal H}^{1}_{X}) 
\to H^{2}(Y_{an},\Z))}{Im(NS(X))}} \rto \dto & 0 \\
0 \rto & 
\displaystyle{\frac{H^{1}(\tilde{F_{an}},\Z)}{Im(H^{1}(\tilde{Y_{an}}
,\Z))^{s}}} \rto &
\displaystyle{\frac{H^{2}(\tilde{Y_{an}},\tilde{F_{an}},\Z)}{Im(NS(X))}} 
\rto & \displaystyle{\frac{Im(H^{2}(\tilde{Y_{an}},\tilde{F_{an}},\Z)
\to H^{2}(\tilde{Y_{an}},\Z))}{Im(NS(X))}} \rto & 0 \enddiagram$$     
The following diagram commutes by functoriality of the extension class 
maps. 

$$\diagram Im(H^{1}(X,{\cal H}_{X}^{1}) \to H^{2}(Y_{an},\Z)) 
\rto^{\;\;\;\psi} \dto &
\displaystyle{J\left(\frac{H^{0}(X_{an},{\cal 
F}_{\sZ})^{s}}{Im(H^{1}(Y_{an},\Z))}\right)} \dto \\ 
Im(H^{2}(\tilde{Y_{an}},\tilde{F_{an}},\Z) \to H^{2}(\tilde{Y_{an}},\Z))
\rto & \displaystyle{J\left(\frac{H^{1} 
(\tilde{F_{an}},\Z)}{Im(H^{1}(\tilde{Y_{an}},\Z))^{s}}\right)} 
\enddiagram$$
where 
$$J\left(\frac{H^{0}(X_{an},{\cal
F}_{\sZ})^{s}}{Im(H^{1}(Y_{an},\Z))^{s}}  
\right) \to   
J\left(\frac{H^{1} 
(\tilde{F_{an}},\Z)}{Im(H^{1}(\tilde{Y_{an}},\Z))^{s}}\right)$$
is induced from the map on the underlying Hodge structures.

Let $$\psi^{'}:Im(H^{1}(X,{\cal H}_{X}^{1}) \to H^{2}(Y_{an},\Z)) \to
J\left(\frac{H^{1}(\tilde{F_{an}},\Z)}{Im(H^{1}(\tilde{Y_{an}},\Z))} 
\right)$$ 
be the composite in the diagram above.

We also have a natural \lq\lq sheaf theoretic\rq\rq map 
$$\phi^{'}: \frac{Im(H^{1}(X,{\cal H}^{1}_{X}) \to
H^{2}(Y_{an},\Z))}{NS(X)} \to 
J\left(\frac{H^{1}(\tilde{F_{an}},\Z)}{Im(H^{1}(\tilde{Y_{an}},\Z))}
\right)$$ 
which is defined as follows. Let $\eta \in Im(H^{1}(X,{\cal H}^{1}_{X})
\to H^{2}(Y_{an},\Z)).$ Consider a line bundle $L$ on $Y$ such that
$c_{1}(L)= \eta.$ Then, $L|_{F}$ (where $F = \cup_{i} F_{i}$) gives an 
element of $Pic^{0}(F)$ and hence an element of $J(H^{1}(F,\Z))$ via the
mapping $Pic^{0}(F) \to J(H^{1}(F,\Z)).$ Under this mapping $NS(X)$
goes to zero, hence we get a well defined mapping
$$\frac{Im(H^{1}(X,{\cal H}^{1}_{X}) \to H^{2}(Y_{an},\Z))}{NS(X)} 
\to J\left(\frac{H^{1}(F_{an},\Z)}{Im(H^{1}(Y_{an},\Z))}\right).$$ 
Now compose with the map (induced by the morphism of the underlying
Hodge structures)
$$J\left(\frac{H^{1}(F_{an},\Z)}{Im(H^{1}(Y_{an},\Z))}\right) \to
J\left(\frac{H^{1}(\tilde{F_{an}},\Z)}{Im(H^{1}(\tilde{Y_{an}}, 
\Z))}\right)$$  
to get
$$\phi^{'}: \frac{Im(H^{1}(X,{\cal H}^{1}_{X}) \to
H^{2}(Y_{an},\Z))}{NS(X)} \to   
J\left(\frac{H^{1}(\tilde{F_{an}},\Z)}{Im(H^{1}(\tilde{Y_{an}},\Z))}
\right).$$

It is easy to see the following diagram commutes
$$\diagram \displaystyle{\frac{Im(H^{1}(X,{\cal H}^1_X) \to 
H^2(Y_{an},\Z))}{Im(NS(X))}} \rto^{\quad\quad\phi} \drto^{\phi^{'}} &
\displaystyle{\frac{H^0(X,{\cal F})}{Im(Pic^0(Y))}} \dto \\
& \displaystyle{J\left(   
\frac{H^1(\tilde{F_{an}},\Z)}{Im(H^1(\tilde{Y_{an}},\Z))}\right)}
\enddiagram$$
where the map 
$\displaystyle{\frac{H^{0}(X,{\cal F})}{Im(Pic^{0}(Y))}} \to 
\displaystyle{J\left( 
\frac{H^{1}(\tilde{F_{an}},\Z)}{Im(H^{1}(\tilde{Y_{an}},\Z))}\right)}$ 
is the composition 
$$\frac{H^{0}(X,{\cal F})}{Im(Pic^{0}(Y))} \to
\frac{Pic^{0}(F)}{Im(Pic^{0}(Y))} \to 
\frac{Pic^{0}(\tilde{F})}{Im(Pic^{0}(\tilde{Y}))} \to  
J\left(\frac{H^{1}(\tilde{F_{an}},\Z)}{Im(H^{1}(\tilde{Y_{an}} 
,\Z))}\right).$$

We now note that the map $H^{0}(X_{an},{\cal F}_{\sZ})^{s} \to
H^{1}(F_{an},\Z)$ is an injective map. Since $W_{0}H^{1}(F_{an},\Q) =
Ker(H^{1}(F_{an},\Q) \to H^{1}(\tilde{F_{an}},\Q))$ and
$H^{0}(X_{an},{\cal F}_{\sZ})^{s}$ is pure of weight one, it follows that 
$$(Ker(H^{1}(F_{an},\Q) \to H^{1}(\tilde{F_{an}},\Q)))
\cap Im(H^{0}(X_{an},{\cal F}_{\sZ}) \to H^{1}(F_{an},\Q)) = 0.$$  
Hence the composite $$H^{0}(X_{an},{\cal F}_{\sZ}) \to
H^{1}(F_{an},\Z) \to
H^{1}(\tilde{F_{an}},\Z)$$ has finite kernel, and is hence injective (as 
${\cal F}_{\sZ} \subset R^{1}\pi^{an}_{*}\Z$ is torsion-free). 
It follows that
$$J\left(\frac{H^{0}(X_{an},{\cal
F}_{\sZ})^{s}}{Im(H^{1}(Y_{an},\Z))^{s}}\right) \to 
J\left(\frac{H^{1}(\tilde{F_{an}},\Z)}{Im(H^{1}(\tilde{Y_{an}},\Z))^{s}}
\right), \hspace{2mm}(+++)$$
has finite kernel.
 
We have the following diagram which shows all the maps we have 
constructed so far (the outer border is not yet known to commute).  
$$\diagram
\displaystyle{J\left(\frac{H^0(X_{an},{\cal F}_{\sZ})^s}
{Im(H^{1}(Y_{an},\Z))^s} \right)} \ddto & &
\displaystyle{
\frac{Im(H^1(X,{\cal H}^{1}_{X}) \to H^{2}(Y_{an},\Z))}{Im(NS(X))}}
\llto_{\psi} \dto^{\phi}|<<\ahook \\ 
& \displaystyle{\frac{J(H^{0}(X_{an},{\cal F}_{\sZ})^{s})}
{Im(Pic^{0}(Y))}} \rto|<<\ahook \dto &
\displaystyle{\frac{H^{0}(X,{\cal F})}{Im(Pic^{0}(Y))}}
\dto \\ 
\displaystyle{J\left(\frac{H^1(\tilde{F_{an}},\Z)}
{Im(H^{1}(\tilde{Y_{an}},\Z))^s}\right)} &
\displaystyle{\frac{J(H^1(\tilde{F_{an}},\Z))}{Im(Pic^0(\tilde{Y}))}}
\lto &
\displaystyle{
\frac{Pic^0(\tilde{F})}
{Im(Pic^0(\tilde{Y}))}} \lto \enddiagram$$

We will prove that the following subdiagram commutes 
$$\diagram 
\displaystyle{J\left(\frac{H^0(X_{an},{\cal
F}_{\sZ})^{s}}{Im(H^{1}(Y_{an},\Z))^{s}}\right)} \dto &
\displaystyle{\frac{Im(H^{1}(X,{\cal H}^{1}_{X}) \to
H^{2}(Y_{an},\Z))}{Im(NS(X))}} \lto_{\psi\quad\quad} \dto^\phi|<<\ahook \\
\displaystyle{J\left(\frac{H^{1}(\tilde{F_{an}},\Z)}
{Im(H^{1}(\tilde{Y_{an}},\Z))^s}\right)} &
\displaystyle{\frac{H^0(X,{\cal F})}{Im(Pic^{0}(Y))}} \lto
\enddiagram(*)$$ 
Note that the composite map 
$$\frac{Im(H^{1}(X,{\cal H}^{1}_{X}) \to H^{2}(Y_{an},\Z)}{Im(NS(Y))}
\to \frac{H^{0}(X,{\cal F})}{Im(Pic^{0}(Y))} \to  
J\left(\frac{H^{1}(\tilde{F_{an}},\Z)}{Im(H^{1}(\tilde{Y_{an}},\Z))}\right)  
$$ is the previously defined
map $\phi^{'},$ and the map 
$$\frac{Im(H^{1}(X,{\cal H}^{1}_{X}) \to H^{2}(Y_{an},\Z)}{Im(NS(Y))}
\to J\left(\frac{H^{0}(X_{an},{\cal 
F}_{\sZ})^{s}}{Im(H^{1}(Y_{an},\Z))^{s}}\right) \to
J\left(\frac{H^{1}(\tilde{F_{an}},\Z)} 
{Im(H^{1}(\tilde{Y_{an}},\Z))}\right)$$  
is our previously defined map $\psi^{'}.$ 
Thus the commutativity of the above diagram is equivalent to proving
$$\psi^{'} = \phi^{'}.$$

Assuming this diagram commutes we finish the proof of Theorem~\ref{(1,1)}
as follows. We claim that the composite map 
$$\frac{Im(H^{1}(X,{\cal H}^{1}_{X}) \to
H^{2}(Y_{an},\Z))}{Im(NS(X))} \to J\left(\frac{H^{1}(\tilde{F_{an}},\Z)}
{Im(H^{1}(\tilde{Y_{an}},\Z))^{s}}\right)$$ 
has finite kernel, as $\phi$ is injective and the map 
$$\frac{H^{0}(X,{\cal F})}{Im(Pic^{0}(Y))} \to
J\left(\frac{H^{1}(\tilde{F_{an}},\Z)}
{Im(H^{1}(\tilde{Y_{an}},\Z))^s}\right)$$ 
has finite kernel (combining lemma~\ref{compare} with the fact that 
the map in $(+++)$ above has finite kernel). 
Thus, $\psi$ has finite kernel. 
Now recall the map $$\psi_{1}:
\Z^{r} \cong \frac{A}{A_{tors}} \to J(M),$$ 
where $A = \displaystyle{\frac{Im(H^{1}(X,{\cal H}^{1}_{X}) \to
H^{2}(Y_{an},\Z))}{Im(NS(X))}}.$ It follows immediately that $\psi_{1}$ has
finite kernel. Since $\displaystyle{\frac{A}{A_{tors}}}$ is a free abelian group, it
follows that $\psi_{1}$ is injective. This is equivalent to proving our
main result, Theorem~\ref{(1,1)}, as has been remarked before.

We now finish the final part of the proof by showing the commutativity
of the diagram$~(*).$ Let $Z$ be a smooth projective variety over $\C$
and $W \subset Z$ be a smooth subvariety. Let $\eta \in H^{2}(Z_{an},\Z)$
be an algebraic class (i.e., let $\eta \in NS(Z)$), such that $\eta
\mapsto 0 \in H^{2}(W_{an},\Z).$ Then, $\eta$ gives rise to the following
pullback diagram 
$$\hspace{-2cm} \diagram 0 \rto &
\displaystyle{\frac{H^{1}(W_{an},\Z)}{H^{1}(Z_{an},\Z)}}
\rto & H^{2}(Z_{an},W_{an},\Z) \rto & Ker(H^{2}(Z_{an},\Z) \to
H^{2}(W_{an},\Z)) \rto & 0 \\ 
0 \rto & \displaystyle{\frac{H^{1}(W_{an},\Z)}{H^{1}(Z_{an},\Z)}} \rto
\udouble & B \rto \uto & \Z[\eta] \rto \uto & 0 \enddiagram$$ 
Thus we have an extension class map 
$$\diagram \Z \cong \Z[\eta] \rto^{\psi^{''}\quad\quad} &
\displaystyle{J\left(\frac{H^{1}(W_{an},\Z)}{H^{1}(Z_{an},\Z)}\right)}.
\enddiagram$$ 
Again given $\eta$ as above consider $L_{\eta},$ a line bundle on $Z$
which has Chern class $\eta.$ Restrict this line bundle on $W$ to get
$L_{\eta}|_{W} \in Pic^{0}(W_{an}) \cong J(H^{1}(W_{an},\Z)).$ This gives
us a well-defined mapping 
$$\diagram \Z \cong \Z[\eta] \rto^{\phi^{''}\quad\quad} &
\displaystyle{J\left(\frac{H^{1}(W_{an},\Z)}{H^{1}(Z_{an},\Z)}\right)}. 
\enddiagram$$    
We now have the following lemma. 

\begin{lemma} \label{abstract} With the above notation $\psi^{''} =
\phi{''},$ i.e., the extension class map and the restriction map
corresponding to the class
$\eta$ are the same.
\end{lemma}

\begin{proof} Consider the following diagram with exact rows and columns 

$$\hspace{-1cm}\diagram 
\rto & H^{0}(W_{an},\Z) \rto \dto & H^{1}(Z_{an},W_{an},\Z) \rto \dto &
H^{1}(Z_{an},\Z) \rto \dto & H^{1}(W_{an},\Z) \rto \dto & \\ 
\rto & \displaystyle{\frac{H^{0}(W_{an},\C)}{F^{1}H^{0}(W_{an},\C)}}
\rto \dto & \displaystyle{\frac{H^{1}(Z_{an},W_{an},\C)}{F^{1} 
H^{1}(Z_{an},W_{an},\C)}} \rto \dto &
\displaystyle{\frac{H^{1}(Z_{an},\C)}{F^{1}H^{1}(Z_{an},\C)}} \rto
\dto &
\displaystyle{\frac{H^{1}(W_{an},\C)}{F^{1}H^{1}(W_{an},\C)}} \rto
\dto & \\ 
\rto & H^{0}(W_{an},{\cal O}_{W_{an}}^{*}) \rto \dto & Pic(Z_{an},W_{an})  
\rto \dto & Pic(Z_{an}) \rto \dto & Pic(W_{an}) \rto \dto & \\ 
\rto & H^{1}(W_{an},\Z) \rto \dto & H^{2}(Z_{an},W_{an},\Z) \rto \dto &
H^{2}(Z_{an},\Z) \rto \dto & H^{2}(W_{an},\Z) \rto \dto & \\ 
\rto & \displaystyle{\frac{H^{1}(W_{an},\C)}{F^{1}H^{1}(W_{an},\C)}}
\rto \dto & \displaystyle{\frac{H^{2}(Z_{an},W_{an},\C)}{F^{1} 
H^{2}(Z_{an},W_{an},\C)}} \rto \dto &
\displaystyle{\frac{H^{2}(Z_{an},\C)}{F^{1}H^{2}(Z_{an},\C)}} \rto
\dto & 
\displaystyle{\frac{H^{2}(W_{an},\C)}{F^{1}H^{2}(W_{an},\C)}} \rto
\dto & \\
& & & & \enddiagram$$  

The above diagram comes from the following $9$-diagram in the category
of sheaves (where the bottom row defines ${\cal O}_{Z_{an}}(-W_{an})^{*},$ 
and $j:Z_{an}-W_{an} \hookrightarrow Z_{an}$ is the inclusion).
$$\diagram  & 0  \dto & 0 \dto & 0 \dto \\
0 \rto & j_{!}\Z_{Z_{an}-W_{an}} \rto \dto & \Z_{Z_{an}} \rto \dto &
\Z_{W_{an}} \rto \dto & 0 \\ 
0 \rto & {\cal O}_{Z_{an}}(-W_{an}) \rto \dto & {\cal O}_{Z_{an}} \rto
\dto & {\cal O}_{W_{an}} \rto \dto & 0 \\  
0 \rto & {\cal O}_{Z_{an}}(-W_{an})^{*} \rto \dto & {\cal O}_{Z_{an}}^{*}  
\rto \dto & {\cal O}_{W_{an}}^{*} \rto \dto & 0 \\
& 0 & 0 & 0 \enddiagram$$

Let $\eta \in H^{2}(Z_{an},\Z)$ such that $\eta \mapsto 0$ both in
$H^{2}(W_{an},\Z)$ and  
$\displaystyle{\frac{H^{2}(Z_{an},\C)}{F^{1}H^{2}(Z_{an},\C)}}$ (i.e.,  
$\eta \in NS(Z)$). By a diagram chase as before we get elements
$\delta_{1}$ and $\delta_{2}$ in the group
$\displaystyle{\frac{H^{1}(W_{an},\C)}{F^{1}H^{1}(W_{an},\C)}}.$ Both
$\delta_{1}$ and $\delta_{2}$ are well-defined in the quotient group 
$\displaystyle{\frac{H^{1}(W_{an},\C)}{F^{1} 
H^{1}(W_{an},\C)}/\left(Im(H^{1}(W_{an},\Z)
+ Im\left(\frac{H^{1}(Z_{an},\C)}{F^{1}
H^{1}(Z_{an},\C)}\right)\right)}.$ 

Now note that  
$$\frac{H^{1}(W_{an},\C)}{F^{1}
H^{1}(W_{an},\C)}/\left(Im(H^{1}(W_{an},\Z) +
Im\left(\frac{H^{1}(Z_{an},\C)}{F^{1}H^{1}(Z_{an},\C)}\right)\right)
\cong J\left(\frac{H^{1}(W_{an},\Z)}{H^{1}(Z_{an},\Z)}\right).$$ Hence
we get 2 maps $\eta \mapsto \bar{\delta_{1}}$ and $\eta \mapsto
\bar{\delta_{2}}$ from  
$$\Z \cong \Z[\eta] \to 
J\left(\frac{H^{1}(W_{an},\Z)}{H^{1}(Z_{an},\Z)}\right)$$ 
where $\bar{\delta}$ denotes the class of $\delta$ in the quotient  
$\displaystyle{J\left( 
\frac{H^{1}(W_{an},\Z)}{H^{1}(Z_{an},\Z)}\right)}.$

We claim that these two maps are nothing but our previously defined
maps ${\phi}^{''}$ and ${\psi}^{''}$ respectively. It is clear that the 
map $\eta \mapsto \bar{\delta_{2}}$ is equal to $\phi^{''}(\eta).$ This is
because we got $\bar{\delta_{2}}$ by first taking a lift of $\eta,$ say 
$\beta_{2}$ in $Pic(Z_{an}),$ then restricting to $Pic(W_{an})$ to get 
$\gamma_{2}$ and finally taking the class $\bar{\delta_{2}} \in 
\displaystyle{J\left(\frac{H^{1}(W_{an},\Z)}{H^{1}(Z_{an},\Z)}\right)}.$
This is exactly how $\phi^{''}(\eta)$ was defined, so $\bar{\delta_{2}}
= \phi^{''}(\eta).$ It is also clear that $\bar{\delta_{1}} = 
\psi^{''}(\eta)$ as the extension class map is defined exactly the same 
way as the map $\eta \mapsto \bar{\delta_{1}}.$  Now by 
lemma~\ref{[PS]}, we have that $\bar{\delta_{1}} = \bar{\delta_{2}}.$  
This implies that $\phi^{''} = \psi^{''}.$ 
\end{proof}

\begin{rmk} In a similar vein, using lemma~\ref{[PS]}, one can show that
the cycle class map with values in Deligne-Beilinson cohomology restricts
to the Abel-Jacobi mapping, on cycles which are homologically trivial.
This is essentially the argument given in [EV], though the need to appeal
to lemma~\ref{[PS]} is not brought out explicitly there. 
\end{rmk}

Let $Z= \tilde{Y}$ and $W= \tilde{F}$ in lemma~\ref{abstract}.  
Let $\eta \in Im(H^{1}(X,{\cal H}^{1}_{X}) \to H^{2}(Y_{an},\Z)),$ 
and let $\eta \mapsto \tilde{\eta} \in H^{2}(\tilde{Y_{an}},\Z).$
Then, $\tilde{\eta}$ satifies the conditions of lemma~\ref{abstract},
i.e., $\tilde{\eta} \in NS(\tilde{Y})$ and $\tilde{\eta} \mapsto 0 \in 
H^{2}(\tilde{F_{an}},\Z).$ Clearly, 
$$\psi^{'}(\eta) = \psi^{''}(\tilde{\eta}),$$ 
and
$$\phi^{'}(\eta) = \phi^{''}(\tilde{\eta}).$$ 
Hence,
$$\psi^{'} = \phi^{'}$$ 
which proves the commutativity of diagram$~(*).$
This finishes the proof of our main result, Theorem~\ref{(1,1)}.

\section{An example}
In this section we give an example of an integral projective variety $X$
over $\C$ which is not normal, for which we have a strict inclusion
\[NS(X)\propsubset \{\alpha\in H^2(X,\Z)\mid\mbox{ $\alpha$ is Zariski
locally trivial and $\alpha_{\sC}\in F^1H^2(X,\Z)$}\}.\]
Our variety will have the property that its normalization $Y$ is
non-singular, and the normalization map $\pi:Y\to X$ is bijective. Then
$H^2(X,\Z)\cong H^2(Y,\Z)$ as mixed Hodge structures, and the subspaces of
Zariski locally trivial classes correspond. Hence the desired property of
$X$ is equivalent to the strictness of the first inclusion
\[NS(X)\propsubset NS(Y)\subset H^2(Y,\Z).\]

We will make use of a variant of a construction in [Ha], III, Ex. 5.9 (see
also [Ha], II, Ex. 5.16b). If $V$ is a non-singular variety over $\C$,
then following [Ha], an {\it infinitesimal extension} of $V$ by an
invertible $\cO_V$-module $\ccL$ is a scheme $W$ with $W_{red}=V$, such
that the nilradical $\cI$ of $\cO_W$ has square zero (so that it is an
$\cO_V$-module), and there is an $\cO_V$-isomorphism $\cI\cong \ccL$. In
other words, there is an exact sequence of sheaves of $\cO_W$-modules
\[0\to \ccL\to \cO_W\to \cO_V\to 0.\]
There is a corresponding exact sequence of sheaves
\[0\to \ccL\to \cO_W^*\to\cO_V^*\to 0,\]
where $\ccL$ is identified with the (multiplicative) subsheaf of units on
$W$ which restrict to $1$ on $V$ (the identification is given on sections
by $s\mapsto 1+s$).  The latter exact sheaf sequence gives rise to an
exact sequence on cohomology
\[ H^1(V,\ccL)\to \Pic W\to \Pic V \by{\delta}  H^2(V,\ccL).\]

The following is an elaboration of [Ha], III, Ex. 5.9 (the proof is left
as an exercise!).
\begin{lemma}\label{infinitesimal}
\begin{points}
\item There is a natural bijection between isomorphism classes of
infinitesimal extensions of $V$ by $\ccL$ and elements
\[\alpha\in H^1(V,\shom_V(\Omega^1_{V/\sC},\ccL)).\]
\item If $W$ is the infinitesimal extension corresponding to $\alpha$, the
boundary map $\delta=\delta_{\alpha}:\Pic V\to H^2(V,\ccL)$ is expressible
as the composition  
\[\Pic V=H^1(V,\cO_V^*)\by{\dlog} H^1(V,\Omega^1_{V/\sC})\by{\cup \alpha} 
H^2(V,\ccL).\]
\item Let $f:\ccL\to\cM$ be a morphism of $\cO_X$-modules, and
$\alpha\mapsto f_*(\alpha)$ under the natural map 
\[f_*:H^1(V,\shom_V(\Omega^1_{V/\sC},\ccL))\to
H^1(V,\shom_V(\Omega^1_{V/\sC},\cM)).\] 
Let $Z$ be the infinitesimal extension of $V$ by $\cM$ determined by
$f_*(\alpha)$. Then there is a unique morphism of $\C$-schemes
$\tilde{f}:Z\to W$, such that the corresponding morphism of reduced
schemes is the identity on $V$, and such that there are commutative
diagrams with exact rows
\[\begin{array}{ccc}
0\to \ccL\to & \cO_W & \to\cO_V\to 0\\
f\downarrow\quad & \tilde{f}^*\downarrow\quad & \veq\quad\\
0\to\cM\to & \cO_Z & \to\cO_V\to 0
\end{array}\]
and
\[\begin{array}{cccr}
H^1(V,\ccL) & \to \Pic W \to & \Pic V & \longby{\delta_{\alpha}}
H^2(V,\ccL)\\  
\quad f_*\downarrow & \tilde{f}^*\downarrow\quad & \veq & \downarrow f_*\quad\\
H^1(V,\cM) & \to \Pic Z \to & \Pic V & \longby{\delta_{f_*(\alpha)}}
H^2(V,\cM)  
\end{array}\]
\end{points}
\end{lemma}

\begin{ex} In the above lemma, take $V=\P^1_{\sC}\times\P^1_{\sC}$,
and $\ccL=\omega_{V}$, $\cM=\cO_{V}$, $f:\ccL\into\cM$ any non-zero
map (since $\omega_{V}\cong\cO_{\sP^1}(-2)\boxtimes\cO_{\sP^1}(-2)$, such
maps $f$ exist). Infinitesimal extensions of $V$ by $\ccL$ are classified
by elements of 
\[H^1(V,\shom_V(\Omega^1_{V/\sC},\omega_V))\cong H^1(V,\Omega^1_{V/\sC}),\]
where we have identified $\shom_V(\Omega^1_{V/\sC},\omega_V)$ with
$\Omega^1_{V/\sC}$ using the non-degenerate bilinear form 
\[\Omega^1_{V/\sC}\tensor_{\cO_V}\Omega^1_{V/\sC}\to\omega_V\]
given by the exterior product of 1-forms. Thus if $\alpha\in
H^1(V,\Omega^1_{V/\sC})$, the corresponding cup-product map
\[H^1(V,\Omega^1_{V/\sC})\longby{\cup\alpha} H^2(V,\omega_V)\]
is just the product with $\alpha$ in the commutative graded ring
\[\oplus_{n\geq 0} H^n(V,\Omega^n_{V/\sC}).\]
Thus if $\alpha$ is the cohomology class of a divisor $D$ on $V$, then for
any divisor $E$ on $V$, we see that by (ii) of the lemma, 
\[\delta_{\alpha}(E)=(D\cdot E)\in \C= H^2(V,\omega_V)\]
is the intersection product of $D$ and $E$ on the non-singular projective
surface $V$. 

We will choose $\alpha$ to be the cohomology class of $D=L_1-L_2$, where
$L_1=\P^1_{\sC}\times\{0\}$ and $L_2=\{0\}\times\P^1_{\sC}$ are elements
of the two rulings on $V=\P^1_{\sC}\times\P^1_{\sC}$.  Since $D^2=-2$,
$\alpha\neq 0$ and $W=(V,\cO(\alpha))$ is a non-trivial infinitesimal
extension of $V$ by $\omega_V$. Note that $H^1(V,\omega_V)=0$, so that
there is an exact sequence
\[0\to \Pic W\to \Pic V \longby{(D\cdot \;)} \Z\to 0\]
($(D\cdot L_2)=1$, so the map to $\Z$ is surjective). Here $\Pic V=\Pic
(\P^1_{\sC}\times\P^1_{\sC})=\Z[L_1]\oplus \Z[L_2]$ is free abelian of rank
2, and as usual, we denote a representative of the class of
$a[L_1]+b[L_2]$ by $\cO_V(a,b)$. 

Next, note that 
\[f_*(\alpha)\in H^1(V,\shom_V(\Omega^1_{V/\sC},\cO_V))=0,\]
since $\Omega^1_{V/\sC}\cong \cO_V(-2,0)\oplus \cO_V(0,-2)$. Hence the
infinitesimal extension $Z$ of $V$ by $\cO_V$ determined by $f_*(\alpha)$
is the trivial extension $(V,\cO_V[\epsilon])$, where $\cO_V[\epsilon]$ is
the sheaf of dual numbers over $\cO_V$. 

There is an obvious way in which we may regard $Z=(V,\cO_V[\epsilon])$ as
a closed subscheme of
$Y=\P^1_{\sC}\times\P^1_{\sC}\times\P^1_{\sC}=V\times\P^1_{\sC}$, whose
underlying reduced scheme is $V\times\{0\}$. 

Finally, we define $X$ to be the $\C$-scheme which is the pushout of $Y$ and
$W$ along the morphisms $\tilde{f}:Z\to W$ and the above inclusion $i:Z\into
Y$, so that there is a commutative pushout diagram
\[ \begin{array}{ccc}
Z & \longby{i} & Y=V\times\P^1_{\sC}\\
\tilde{f} \downarrow &&\quad \downarrow \pi\\
W & \longby{j} & X
\end{array}\]
Since $\tilde{f}$ is a finite and bijective morphism, we see that for each
affine open $U=\Spec A$ in $Y$, $U\cap Z=\Spec A/I$ is affine, and finite
over the affine open subscheme $\tilde{f}(U\cap Z)=\Spec B\subset W$. The
image of $U$ in $X$ is then defined as the affine scheme $\Spec C$, where
$C$ is the inverse image of $B$ in $A$ under the surjection $A\onto A/I$.
One shows easily that $C$ is in fact a finitely generated $\C$-subalgebra of
$A$, and $A$ is a finite $C$-module with conductor ideal $I$. Further, the
construction of $C$ localizes well. Hence the local schemes $\Spec C$ can
be glued together to yield the $\C$-scheme $X$. 

We claim that this scheme $X$ has the desired properties, \ie $X$ is an
integral projective scheme over $\C$ with normalization $\pi:Y\to X$, such
that  $Y$ is non-singular and bijective with $X$, while $NS(X)\to NS(Y)$
is a strict inclusion.  

That $X$ is integral with $Y$ as its normalization is clear, from the
description of its affine open sets above. Next, since $\alpha=[D]$, and
$(D\cdot (L_1+L_2))=0$, there is a unique $\ccH\in\Pic W$ such that
$\ccH\tensor\cO_V=\cO_V(1,1)$. Also $\Pic Z\to\Pic V$ is an isomorphism.
We have $\Pic Y=\Z^{\oplus 3}$, where we may regard the restriction map
$\Pic Y\to\Pic V\times\{0\}=\Pic V=\Z^{\oplus 2}$ as projection on the
first 2 factors. Hence we see that the very ample invertible sheaf
$\cO_Y(1,1,1)$ has the property that there is an isomorphism
$\tilde{f}^*\ccH\cong \cO_Y(1,1,1)\tensor\cO_Z$. From the Mayer-Vietoris
sequence of sheaves of rings 
\[0\to \cO_X\to \pi_*\cO_Y\oplus j_*\cO_W \to (i\circ \pi)_*\cO_Z\to 0\]
we have a corresponding sequence of sheaves of unit groups
 \[0\to \cO_X^*\to \pi_*\cO_Y^*\oplus j_*\cO_W^* \to (i\circ
\pi)_*\cO_Z^*\to 0\] 
leading to an exact sequence
\[H^0(Y,\cO_Y^*)\oplus H^0(W,\cO_W^*)\to H^0(Z,\cO_Z^*)\to \Pic X \to \Pic
Y\oplus \Pic W \to \Pic Z.\] 
Hence there exists an invertible sheaf $\cA$ on $X$ with
$\pi^*\cA=\cO_Y(1,1,1)$ (and $j^*\cA=\ccH$). Since $\pi$ is  finite, and
$\pi^*\cA$ is ample on $Y$, we have that $\cA$ is ample on $X$. Hence $X$
is projective.  

>From the exact sequence
\[0\to\omega_V\to\cO_W\to \cO_V\to 0,\]
we see that $H^0(W,\cO_W)\to H^0(V,\cO_V)=\C$ is an isomorphism. On the
other hand, we see at once that $H^0(Z,\cO_Z)=\C[\epsilon]$ is the ring of
dual numbers. Hence we get analogous formulas for the unit groups. Thus
there is an exact sequence 
\[0\to \C\to \Pic X\to \Pic Y\oplus \Pic W \to \Pic Z \to 0\] 
(note that $\Pic Z=\Pic V$ is a quotient of $\Pic Y$).
Since $\Pic W\into \Pic Z$ is an inclusion $\Z\into \Z^{\oplus 2}$ as a
direct summand, while $\Pic Y\to\Pic Z$ is the projection $\Z^{\oplus
3}\onto \Z^{\oplus 2}$, we see that ${\rm image}\,(\Pic X\to\Pic Y=NS(Y))$
is a direct summand $\Z^{\oplus 2}\into \Z^{\oplus 3}$. Thus
$NS(X)=\Z^{\oplus 2}$ is strictly contained in $NS(Y)=\Z^{\oplus 3}$. 
$\square$
\end{ex}

\end{document}